\documentclass[11pt,twoside]{article}
\usepackage[dvips]{graphics}
\usepackage{amsmath}
\usepackage{amsmath,graphicx,amsfonts,amsrefs,amssymb}
\usepackage{setspace, enumitem, mathtools, float, yfonts}
\usepackage{graphicx, algorithm, algpseudocode}
\usepackage[table]{xcolor}
\newtheorem{theorem}{Theorem}[section]

\def\qed{\hfill\rule{1ex}{1ex}\\}

\begin{document}
	
	\title{Spectral methods for prescribed mean curvature equations}
	\author{Jonas Haug\footnote{Department of Mathematics, Texas State University, qcd5@txstate.edu}, Rachel Jewell\footnote{Department of Computer Science, Texas State University, kju11@txstate.edu}, and  Ray Treinen\footnote{Department of Mathematics, Texas State University, 601 University Dr., San Marcos, TX 78666, rt30@txstate.edu}
	}
	\maketitle
	
	\begin{abstract}
		We consider  prescribed mean curvature equations whose solutions are minimal surfaces, constant mean curvature surfaces, or capillary surfaces.  We consider both Dirichlet boundary conditions for Plateau problems and nonlinear Neumann boundary conditions for capillary problems and we consider domains in $\mathbf{R}^2$ to be rectangles, disks, or annuli.
		
		We present spectral methods for approximating  solutions of the associated boundary value problems.  These are either based on Chebyshev or Chebyshev-Fourier methods depending on the geometry of the domain.  The non-linearity in the prescribed mean curvature equations is treated with a Newton method.  The algorithms are designed to be adaptive; if the prescribed tolerances are not met then the resolution of the solution is increased until the tolerances are achieved. 
	\end{abstract}
	
%
	
	
	\section{Introduction}
	\label{intro}
	
	We develop computational methods for approximating solutions to some prescribed mean-curvature equations over certain domains in $\mathbb{R}^2$.  We will consider  domains $\Omega\subset \mathbb{R}^2$ to be rectangles, disks, or annuli.  For our problems,
	$u:\Omega\rightarrow \mathbb{R}$ is interpreted as a height over the base domain $\Omega$, where $u$ satisfies
	\begin{equation}\label{eqn:pde}
		\nabla \cdot \frac{\nabla u}{\sqrt{1 + |\nabla u|^2}} = \mathcal{H}(u) \mbox{ in } \Omega,
	\end{equation}
	where the left side of the equation is the mean curvature operator.
	Here $\mathcal{H}(u) = 0$  corresponds to the minimal surface equation, $\mathcal{H}(u) = 2H\in\mathbb{R}$ fixed corresponds to a constant mean curvature (CMC) equation, and $\mathcal{H}(u) = \kappa u$ corresponds to a capillary surface with capillary constant $\kappa >0$.   Physically capillary surfaces are given by the height of the interface between two fluids in equilibrium.  Here $\kappa = \rho g/\sigma$, with the density difference $\rho$ measured between the two fluids, $g$ is the gravitational constant, and $\sigma$ is the surface tension for that fluid-fluid interface.  Traditionally, scaling arguments are used to simplify the mathematical problem so that $\kappa = 1$.  In what follows we keep the general parameter available, but we use $\kappa = 1$ in our numerical results.
	
	Generically, we are interested in boundary conditions that are either the Dirichlet data
	\begin{equation}
		u = g \mbox{ on } \partial \Omega
	\end{equation}
	for some prescribed function $g$, or the so-called capillary data
	\begin{equation}\label{eqn:capdata}
		\frac{\nabla u}{\sqrt{1 + |\nabla u|^2}} \cdot \mathbf{n} = \cos\gamma \mbox{ on } \partial \Omega
	\end{equation}
	where $\mathbf{n}$ is the outward unit normal to $\Omega$, and $\gamma\in[0,\pi]$ is the contact angle as measured from within a fluid at the contact with the wall of the tube $\partial\Omega \times \mathbb{R}$.    While one can think of the second type of boundary data as nonlinear Neumann data, it is also known as the natural boundary condition for a capillary surface, as first derived by Gauss in 1830 \cite{Gauss}.  It is, of course, possible to mix these types of boundary data as well.
	
This paper is organized as follows.	The existence and uniqueness theory for the different combinations of equations and boundary data is not straightforward, with restrictions on both the type and quantities in the boundary data as well as the geometry of $\Omega$.  We will  survey these results in Section~\ref{Lit}.  In Section~\ref{operators} we provide details for the operators associated with the various combinations of PDEs, boundary conditions, and geometries of the domains.  This includes the Fr\'echet derivatives of the operators to yield an iterative linearization of these operators suitable for using Newton's method to solve these nonlinear boundary value problems.	In Section~\ref{spectral} we discuss the Chebyshev and Chebyshev-Fourier spectral methods for the  linearized problems at each iterative step of Newton's method.	With all of this developed in a modular fashion, in Section~\ref{newton} we are able to present a unified approach to our adaptive Newton's method.  Finally, in Section~\ref{results} we provide some examples of our algorithms  highlighting the strengths and limitations of this approach.

We have released Matlab implementations of these algorithms under an open source license.  These files are available in a GitHub repository at \\
	\texttt{https://github.com/raytreinen/Spectral-Methods-Mean-Curvature}\\

\section{A discussion of appropriate boundary value problems}
\label{Lit}
	
\begin{figure}[t!]
		\centering
		\scalebox{0.195}{\includegraphics{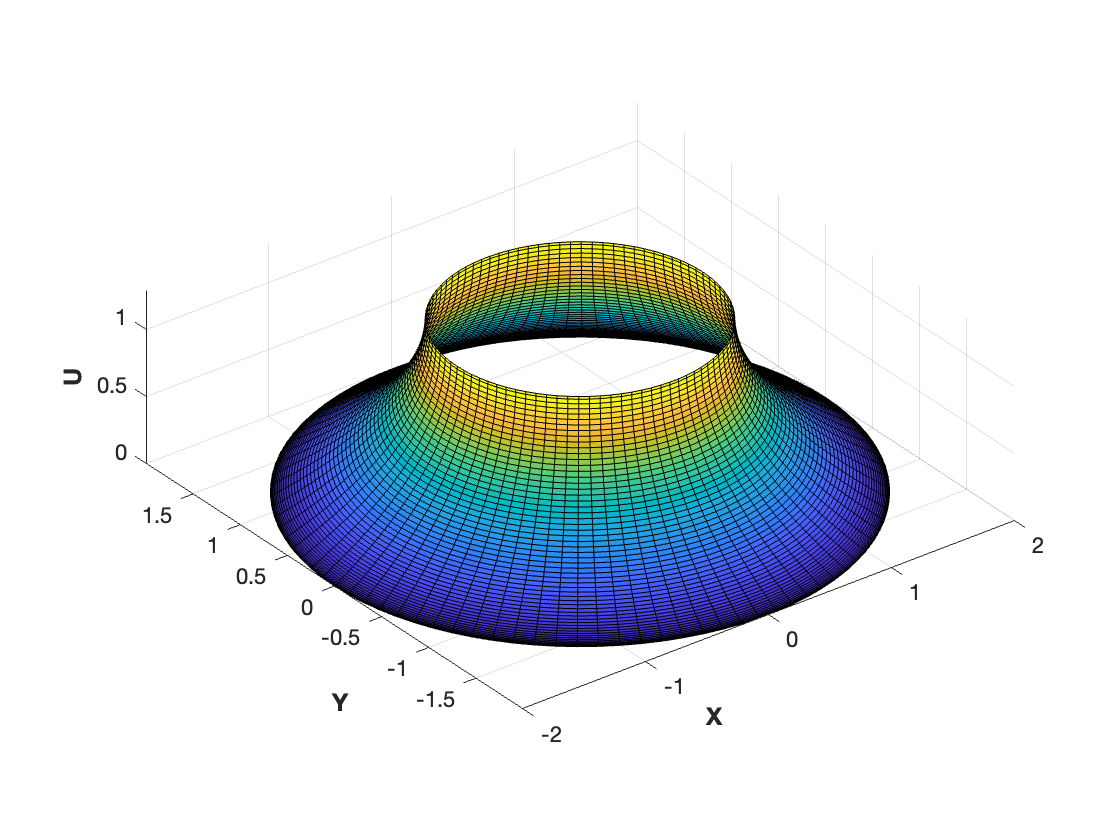}}
		\scalebox{0.195}{\includegraphics{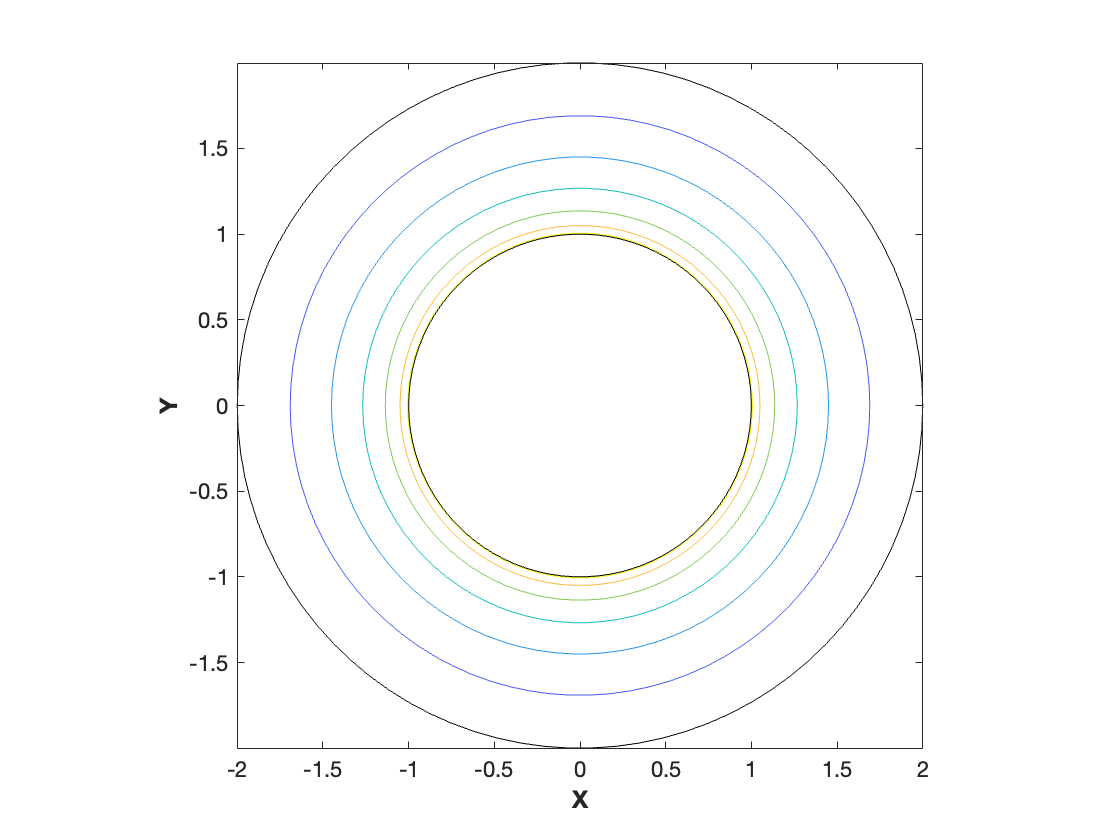}}
		\caption{A minimal surface over an annular domain with a gradient close to blow-up.  The left is a plot of the surface height and the right is a contour plot.  This provides an illustration of the possible restrictions on the arbitrary boundary data for solutions to exist for non-convex domains.  The height of the surface is prescribed to be 0 at the outer radius $b = 2$ and is prescribed to be 1.28792 at the inner radius $a = 1$.  If much larger heights are prescribed at $a$ the solution is no longer described by a graph over the base domain, but can be seen to be some catenoid for a further range of heights.  While not the only treatment of these results, Bliss \cite{Bliss1925} describes the process and the resulting stability analysis nicely. }
		\label{fig:MinDirAnnCat}
	\end{figure}

	The Plateau problems are traditionally considered when \eqref{eqn:pde} is coupled with  Dirichlet boundary data $u = g$  on $\partial \Omega$ for $\mathcal{H}(u) = 0$ or $\mathcal{H}(u) = 2H$.  	The introduction of L\'opez's book \cite{Lopez2013} provides a nice literature review, and we sample from that work.  There are restrictions on the values of $H$ for Plateau problems for CMC equations, and the divergence theorem implies  
	$$
	|2H| \leq \frac{|\partial \Omega| }{|\Omega|}.
	$$
	Finn \cite{Finn1965} showed
	\begin{theorem}
	For the Plateau problem with $H=0$, there exists a minimal graph on $\Omega\subset\mathbb{R}^2$ for each arbitrary boundary condition $g$ if and only if $\Omega$ is convex.
	\end{theorem}
We show an example of a Plateau problem in Figure~\ref{fig:MinDirAnnCat}, though the domain is not convex.	For the cases where $H\ne 0$, Serrin \cite{Serrin1969} proved
	\begin{theorem}
	Let $\Omega\subset\mathbb{R}^2$ is a bounded strictly convex domain, and let $\mathcal{H}(u) = 2H \ne 0$.  Then there exists a solution of \eqref{eqn:pde} with arbitrary boundary data $g$ if and only if
	$$
	k \geq 2H > 0,
	$$
	where $k$ is the curvature along the plane curve $\partial \Omega$.
	\end{theorem}
	
	The pairing of $\mathcal{H}(u) = \kappa u$ with Dirichlet boundary data is less common, but can certainly be considered.  However, the so-called capillary boundary data \eqref{eqn:capdata} is often paired with $\mathcal{H}(u) = 2H$ to describe the equilibrium shapes of  fluids in the absence of gravity, or with $\mathcal{H}(u) = \kappa u$ for those fluids in a gravity well.  Finn's monograph \cite{ecs} is the standard reference for these problems.
	
	We first discuss the $\mathcal{H}(u) = 2H$ case for these capillary problems.  Another application of the divergence theorem yields
\begin{equation}\label{eqn:2H}
	2H = \frac{1}{|\Omega|}\int_{\partial\Omega} \cos(\gamma(s))\, ds.
\end{equation}
One typically studies the constant contact angle problem, which leads to 
$$
	2H = \frac{|\partial\Omega|}{|\Omega|}\cos\gamma.
$$
	Either of these completely determines the mean curvature from the geometry of the domain and the data $\cos\gamma$.  This is further explored in \cite{ecs}, but we mention in particular the following two results
	\begin{theorem}
		\label{thm:cmc_corner}
	If $\Omega$ contains a corner with interior angle $2\alpha$, and if $\alpha + \gamma < \pi/2$, then when $\mathcal{H}(u) = 2H$ and \eqref{eqn:2H} holds, the problem \eqref{eqn:pde} with \eqref{eqn:capdata} admits no solution.
		\end{theorem}
		
			\begin{figure}[t]
					\centering
	\scalebox{0.195}{\includegraphics{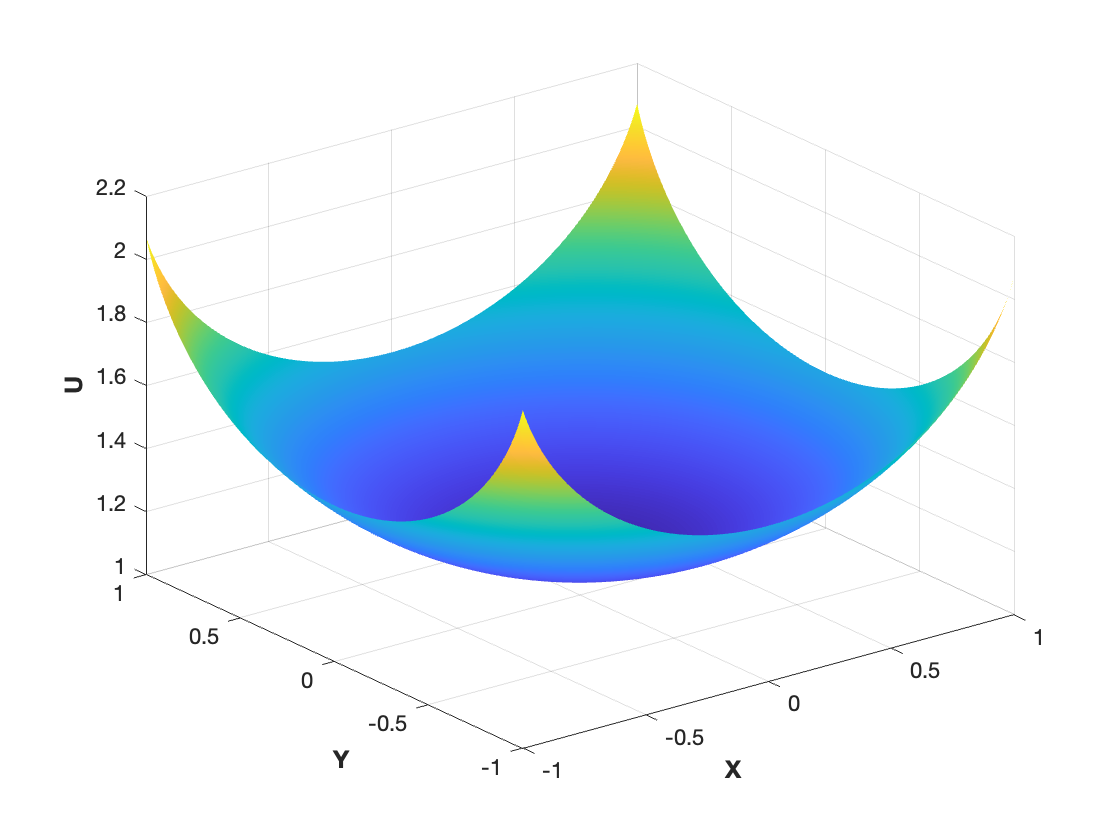}}
\scalebox{0.195}{\includegraphics{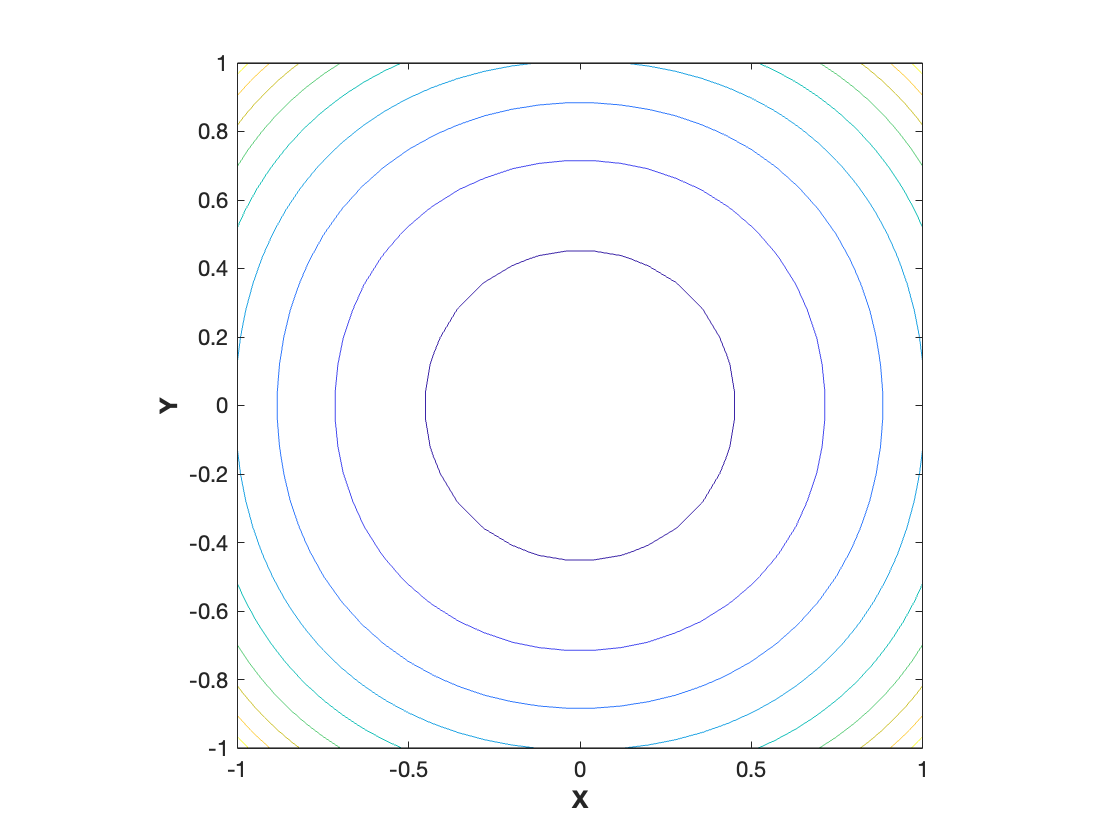}}
					\caption{A capillary surface over a square domain with a near-optimally extreme contact angle of $\gamma = \pi/4 + 0.035$.  The surface height plot on the left shows the gradients getting large, and the contour plot on the right  is a helpful companion figure for the proof of blow-up of capillary surfaces as found in \cite{ecs}.  }
					\label{fig:CapCorner}
				\end{figure}

		\begin{theorem}\label{CMCunstable}
		Capillary surfaces without gravity  are always unstable with respect to boundary perturbations, in the sense that when $\mathcal{H}(u) = 2H$ and \eqref{eqn:2H} holds, for any boundary  $\Sigma = \partial\Omega$ for which a solution of \eqref{eqn:pde}, \eqref{eqn:capdata} exists, there is an arbitrarily small perturbation $\Sigma \rightarrow \hat\Sigma$ such that $\hat\Sigma$ admits no solution.  If $\gamma = 0$, the perturbation can be chosen arbitrarily small, both in $\Sigma$ and in the unit normal vector to $\Sigma$.
		\end{theorem}
		
		For capillary surfaces in a vertical gravity field, where $\mathcal{H}(u) = \kappa u$, Finn \cite{ecs} collects the existence and uniqueness theory that applies to the examples considered in this paper.    The character of Theorem~\ref{thm:cmc_corner} is preserved in the gravity case, and we refer the reader to that work for the technical details.  For our problems with rectangular domains, $\alpha = \pi/4$, so we restrict our contact angles to satisfy $\gamma > \pi/4$ so that  solutions exist, and we note that as $\gamma \rightarrow \pi/4$, the gradient blows up.  This is illustrated in Figure~\ref{fig:CapCorner}.
	
When these boundary conditions are applied to the minimal surface equation, an application of the divergence theorem results in
	\begin{eqnarray}
		0 &=& \int_\Omega \nabla \cdot \frac{\nabla u}{\sqrt{1 + |\nabla u|^2}} \, dA \nonumber\\
		&=& \int_{\partial\Omega} \cos\gamma(s) \, ds. \label{eqn:int}
	\end{eqnarray}
So we see that the total flux along the boundary must be zero.	This condition is the same kind as in Theorem~\ref{CMCunstable}, and we mention Ros \cite{Ros1996} as a place to read more in the setting of surfaces forces applied in the direction of the conormal along the boundary of the surface.

	\section{Operators and Fr\'echet derivatives}
	\label{operators}
	
	In order to implement our method we will need to expand \eqref{eqn:pde} with the different choices for $\mathcal{H}(u)$.  For rectangular domains we use
	\begin{equation}
	\nabla \cdot \frac{\nabla u}{\sqrt{1 + |\nabla u|^2}} = \frac{1}{(1 + |\nabla u|^2)^{3/2}}\left((1 + u^2_y) u_{xx}  - 2u_xu_yu_{xy} +(1 + u^2_x) u_{yy}  \right),
	\end{equation}
	and we define $M$ by
	\begin{equation}
		M(u) = (1 + u^2_y) u_{xx}  - 2u_xu_yu_{xy} +(1 + u^2_x) u_{yy} .
	\end{equation}
	Then the minimal surface equation is $M(u) = 0$.   For the CMC equation we then have
	\begin{equation}
	M(u) -2H (1 + |\nabla u|^2)^{3/2} = 0,
	\end{equation}
	and for capillary surfaces we have
	\begin{equation}
		M(u) - \kappa u (1 + |\nabla u|^2)^{3/2} = 0.
	\end{equation}
	
We will generically refer to these nonlinear operators as $N(u)$, with solutions satisfying $N(u) = 0$.  Later we will include our boundary conditions in this operator.  We will then use the Fr\'echet derivative of each of our choices of $N(u)$, given by $N^\prime = \partial N / \partial u$.  This was discussed by Aurentz and Trefethen \cite{AurentzTrefethen2017} in the setting of lower dimensional problems using spectral methods, and they refer to Zeidler \cite{Zeidler1986}, which is a nice reference for further reading.  We define the differential operators $D_{xx}, D_{yy}, D_{xy}, D_x$, and $D_y$ in the standard way. Then we can use the Fr\'echet derivatives in a modular definition by setting
	\begin{equation}
	F(u) = (1 + u^2_y) D_{xx}  - 2u_xu_yD_{xy} +(1 + u^2_x) D_{yy} + 2(u_{yy}u_x - u_{xy}u_y)D_x + 2(u_{xx}u_y - u_{xy}u_x)D_y.
	\end{equation}
	The Fr\'echet derivative of the nonlinear minimal surface operator is then $F(u)$, where $F$ is evaluated at the function $u$, giving information of the linear behavior of $M$ at that particular $u$ in the sense that $F(u)u$ is a linearization of $M$ at $u$.  We will use this in Newton's method below.
	We first use this to compute the Fr\'echet derivative for the CMC equation modularly as
	\begin{equation}
	F(u) - 6H \sqrt{1 + |\nabla u|^2} \left( u_xD_x + u_yD_y  \right),
	\end{equation}
	and again for the capillary equation as
	\begin{equation}
	F(u) - \kappa (1 + |\nabla u|^2)^{3/2} - 3\kappa u \sqrt{1 + |\nabla u|^2} \left( u_xD_x + u_yD_y  \right).
\end{equation}	

The boundary conditions come in two varieties.  The Dirichlet data
\begin{equation}
u - g = 0 \mbox{ on } \partial \Omega
\end{equation}
has the resulting Fr\'echet derivative of 1.  The nonlinear Neumann data, or capillary data is
\begin{equation}\label{eqn:capbc}
\mathbf{n}\cdot \frac{\nabla u}{\sqrt{1 + |\nabla u|^2}} = \cos\gamma
\end{equation}
on $\partial\Omega$ with outward normal $\mathbf{n} = \langle n_1, n_2 \rangle$, and we will write this as
\begin{equation}
n_1u_x + n_2u_y - \cos\gamma\sqrt{1 + u^2_x + u^2_y} = 0.
\end{equation}
This has the resulting Fr\'echet derivative
\begin{equation}
\left(  n_1 - \frac{\cos\gamma u_x}{\sqrt{1 + u_x^2 + u^2_y}}    \right)D_x + \left(  n_2- \frac{\cos\gamma u_y}{\sqrt{1 + u_x^2 + u^2_y}}    \right)D_y.
\end{equation}
These boundary conditions are well defined for the entire boundary except at the corners.  There we use the measure theoretic normal with $|\mathbf{n}| = 1/\sqrt{2}$.  This is the limit of the average of the exterior normal vectors to $\Omega$ in the following sense.
At each corner intersect $\Omega$ with balls of decreasing radii, each centered at the corner, and denote this set by $\Omega_r = \Omega\cap B_r$.  Then take the average of the normals to $\partial\Omega$ on $\Omega_r$ and take the limit as $r\rightarrow 0$.
See Giusti \cite[p 44]{Giusti1984}.  A recent work by Burns, Fortunato, Julien, and Vasil  \cite{Burns2022} also treats the corner boundary conditions, however, the implementation they derive does not hold for our operators.

Then for problems on a disk or annulus we use cylindrical coordinates with radial variable $r$ and angular variable $\theta$.  We follow the same basic structure as before, and traditionally
\begin{equation}
\tilde M(u) = \left( 1 + \frac{1}{r^2} u^2_\theta  \right) u_{rr} + \frac{1}{r^2}\left( 1 + u^2_r \right) u_{\theta\theta} - \frac{2}{r^2} u_r u_\theta u_{r\theta} + \frac{2}{r^3} u_r u^2_\theta + \frac{1}{r}\left( u_r + u^3_r\right).
\end{equation}
We use $M(u) = r^3\tilde M(u)$, or
\begin{equation}
	M(u) =  r\left( r^2 + u^2_\theta  \right) u_{rr} + r\left( 1 + u^2_r \right) u_{\theta\theta} - 2r u_r u_\theta u_{r\theta} + 2 u_r u^2_\theta + r^2\left( u_r + u^3_r\right).
\end{equation}
which gives us a polar form of the minimal surface equation $M(u) = 0$.  We are reusing the same notation of operator $M$, and $F$ to follow, and the context will be sufficient to differentiate between the rectangular and polar forms of the same equations.  Then the CMC equation becomes
\begin{equation}
M(u) - 2H \left( r^2(1 + u^2_r) + u^2_\theta \right)^{3/2} = 0,
\end{equation}
and for capillary surfaces we have
\begin{equation}
	M(u) - \kappa u \left( r^2(1 + u^2_r) + u^2_\theta \right)^{3/2} = 0.
\end{equation}
Then, as before, we compute the Fr\'echet derivatives.  First we obtain
\begin{eqnarray}
F(u) &=& \left( 2 r u_r u_{\theta\theta}  - 2 r u_\theta u_{r\theta} + 2 u^2_\theta + r^2(1 + 3 u^2_r)  \right) D_r + \left( 2 r u_\theta u_{rr} - 2 r u_r u_{r\theta} + 4 u_r u_\theta  \right) D_\theta \nonumber\\
& & + r \left( r^2 +  u^2_\theta  \right) D_{rr} + r \left( 1 + u^2_r  \right) D_{\theta\theta} - 2r u_r u_\theta D_{r\theta}.
\end{eqnarray}
Then we use this to compute the Fr\'echet derivative for the CMC equation as 
\begin{equation}
F(u) -6H \sqrt{r^2(1 + u^2_r) + u^2_\theta}\left( r^2 u_r D_r +  u_\theta D_\theta \right)
\end{equation}
and then for the capillary equation as
\begin{equation}
F(u) - 3\kappa u \sqrt{r^2(1 + u^2_r) + u^2_\theta}\left( r^2 u_r D_r +  u_\theta D_\theta \right) - \kappa \left( r^2(1 + u^2_r )+ u^2_\theta \right)^{3/2}.
\end{equation}
The equation \eqref{eqn:capbc} becomes less complicated in polar form, as the exterior normal to $\Omega$ is in the radial direction.  Thus the equation for the boundary conditions can be written as
\begin{equation}
\pm r u_r - \cos\gamma\sqrt{r^2 + r^2 u^2_r + u^2_\theta} = 0.
\end{equation}
For the disk or for the outer radius of an annulus we have $+$ in $\pm$, and for the inner radius of an annulus we have $-$ in $\pm$.  
It follows that the Fr\'echet derivative is 
\begin{equation}
 \left(\pm  r - \frac{r^2 u_r \cos\gamma}{\sqrt{r^2 + r^2u^2_r + u^2_\theta}} \right) D_r - \frac{ u_\theta \cos\gamma}{\sqrt{r^2 + r^2u^2_r + u^2_\theta}} D_\theta
\end{equation}
with the same parity of $\pm$ applied to the geometry at hand.

We will be including the boundary conditions in the operator $N(u)$ and also in the Fr\'echet derivatives, and we will discuss the details in the next section.  Once we have done that then we will define the linear operator for the minimal surface boundary value problems  to be $Lu = F(u)u$  with matrix $L$, and we modularly update $L$ to the other problems as appropriate.

\section{Spectral methods}
\label{spectral}

We use Chebyshev and Chebyshev-Fourier spectral methods to solve our equations, and we include a brief overview of these topics here.  For further reading see Trefethen \cite{Trefethen2000}, and the updates found in Chebfun \cite{Chebfun}.

We will first describe differentiation matrices based on functions $f:[-1,1]\rightarrow \mathbb{R}$.  Chebyshev differentiation matrices in this setting are based on polynomial interpolation of sample points of $f$ and $f^\prime$ at Chebyshev grid points $x_j = \cos(\theta_j)\in[-1,1]$ with angles $\theta_j$ equally spaced angles over $[0,\pi]$.  These matrices can be realized as representing the linear transformation between two vectors of data corresponding those grid points points, denoted by $\mathbf{f}$ being mapped to $\mathbf{f^\prime}$.  Fourier differentiation matrices in this setting are based on interpolation by trigonometric polynomials at equally spaced points. They are best used when applied to periodic problems and are also much more commonly studied.

We start with rectangular domains, and, for convenience, we build these $n$  points using Chebfun:
\begin{verbatim}
	x = chebpts(n,[a;b]);
	y = chebpts(n,[c;d]);
\end{verbatim}
then we have two differentiation matrices, $D$ and $D^2$,  implemented by the Chebfun commands
\begin{verbatim}
	D1 = diffmat(n,1,X);
	D2 = diffmat(n,2,X);
\end{verbatim}
where $X = [a,b]$ for derivatives with respect to $x$ and $X = [c,d]$ if we are taking derivatives with respect to $y$, and the input 1 or 2 indicates the number of derivatives.  Then we use Kronecker products to extend these ordinary derivatives to partial derivatives on the rectangle.  Our function height $u$ over the rectangle will be represented as a vector in the usual way of stripping off columns of $n\times n$ data from left to right and concatenating the data into a long $n^2\times 1$ vector: in Matlab $\mathbf{u} = \mathbf{u}(:)$. Then the matrix that gives the discrete partial derivative with respect to $x$ is given by $D\otimes I$ where $I\in \mathbb{R}^{n\times n}$ is the identity matrix and the discrete partial derivative with respect to $y$ is given by $I\otimes D$.  This generalizes to higher derivatives by replacing the differentiation matrix with the higher order version, and of course we implement this in Matlab with the matrices we have generated above using Chebfun.  To get mixed derivatives we simply multiply $D_xD_y$.

As we described in the last section, we also study problems using cylindrical coordinates.  The polar form of Chebyshev grid points has some subtleties when the domain is a disk.  The complication is that if one takes the domain to be $[0,1]\times [-\pi,\pi]$ with (scaled and translated) Chebyshev points $\mathbf{r}$ on $[0,1]$ and evenly spaced points $\mathbf{\theta}$ on $[-\pi,\pi]$, then there are too many points clustered at the origin.  There is also the problem that for the differential equations we are considering the polar form has division by the radial variable and the radius 0 is in each ray from the origin.  To counter these problems we use a double covering of the disk.  This idea is not new, and we refer to Boyd, \cite{Boyd2001}, Fornberg \cite{Fornberg1996},  Trefethen \cite{Trefethen2000} and Wilber, Townsend, and Wright \cite{WilberTownsendWright2017} for further reading.  We will use $(\mathbf{r},\mathbf{\theta}) \in [-1,1]\times [-\pi,\pi]$ with an even number of Chebyshev points $\mathbf{r}$ in $[-1,1]$ so that the value of 0 is skipped.

Wilber, Townsend, and Wright \cite{WilberTownsendWright2017} describe the setting with some detail, and we summarize their discussion of the disk analogue of the double Fourier sphere method.  We generically consider $f(r,\theta)$ defined on $[-1,1]\times [-\pi,\pi]$ and we define $g(r,\theta)$ and $h(r,\theta)$ on $[0,1]\times [-\pi,\pi]$ with $g(r,\theta) = f(r,\theta - \pi)$ and $h(r,\theta) = f(r,\theta)$.  Then the doubled expression for $f$, denoted by $\tilde f$ is given by
\begin{equation}
	\label{eqn:BMC}
\tilde f = \left\{
\begin{array}{ll}
	g(r,\theta + \pi)& (r,\theta) \in [0,1]\times[-\pi,0], \\
	h(r,\theta) & (r,\theta) \in [0,1]\times[0,\pi],  \\
	g(-r,\theta) & (r,\theta) \in [-1,0]\times[0,\pi],  \\
	h(-r,\theta + \pi)) & (r,\theta) \in [-1,0]\times[-\pi,0]. 
\end{array} 
\right.
\end{equation}
Then they observe that $\tilde f$ possesses block-mirror centrosymmetric (BMC) structure and they define functions satisfying \eqref{eqn:BMC} to be BMC functions.  Then to enforce continuity on the disk they further define BMC-II functions to be BMC functions that additionally satisfy that $f(0,\cdot) = \alpha$ for some constant $\alpha$.  With this established, we only need to work with subsets of the doubled data given by $g$ and $h$ to represent the data on the entire disk.

To build differentiation matrices using this double covering we also refer to Chapter~11 of Trefethen's book \cite{Trefethen2000}.  There the Chebyshev points are reversed, as the first entry in a vector $\mathbf{x}$ of those grid points corresponds to $x = 1$, and the last entry corresponds to $x = -1$.  We use Chebfun's more natural ordering of the grid points.  Then in our setting we have the first derivative matrix
\begin{equation}
D_r = 
\begin{bmatrix}
E_1 & E_2 \\
E_3 & E_4	
\end{bmatrix}
\end{equation}
where the blocks line up with the functions in \eqref{eqn:BMC}, and similarly for the second derivative matrix
\begin{equation}
	D_{rr} = 
	\begin{bmatrix}
		D_1 & D_2 \\
		D_3 & D_4	
	\end{bmatrix}
\end{equation}
for the double covered $\tilde f$.  Let  $m$ be the number evenly spaced angular grid points in the double covering, $n$ the number of Chebyshev points, and $I$ the identity matrix of size $m/2$.  Then we extract the halves of the differentiation matrices that correspond to positive values of $r$ by setting 
\begin{equation}
D_r = 
\begin{bmatrix}
0 & I \\
I & 0
\end{bmatrix} 
\otimes E_3 + I_m\otimes E_4
\end{equation}
and 
\begin{equation}
	D_{rr} = 
	\begin{bmatrix}
		0 & I \\
		I & 0
	\end{bmatrix} 
	\otimes D_3 + I_m\otimes D_4.
\end{equation}
We build the Fourier differentiation matrices $D_\theta$ and $D_{\theta\theta}$ and the mixed matrix $D_{r\theta}$ in Matlab using Chebfun by
\begin{verbatim}
D1t = diffmat(m, 1, 'periodic', [-pi,pi]);
D2t = diffmat(m, 2, 'periodic', [-pi,pi]);
\end{verbatim}
and 
\begin{verbatim}
Dth = kron(D1t,eye(n/2));
Dthth = kron(D2t,eye(n/2));
Drth = Dr * Dth;.
\end{verbatim}

The annular domain is, of course, similar to what we have just outlined, but we do not have a double cover.  We use the full Chebyshev differentiation matrix $D$ over the radial interval $[a,b]$ and we have the Fourier differentiation matrix $D_\theta$ over $\theta\in[-\pi,\pi]$.  Then we use Kronecker products $D_r = I_m\otimes D$, $D_{rr} = I_m\otimes D^2$, $D_\theta = D\otimes I_n$, $D_{\theta\theta} = D^2\otimes I_n$, and $D_{r\theta} = D_rD_\theta$.  Here the context indicates if $D$ corresponds to Chebyshev or Fourier differentiation matrices.  As before, we use Chebfun's diffmat to build these component matrices and we explicitly build the second derivative matrices with diffmat.

With the tools we have just described we can assemble the operators $N$ and $L$ for each of the problems we have discussed.  

To implement the boundary conditions, we keep track of the appropriate boundary elements of our height vector $\mathbf{v}$ and for the rows of our vector $N(\mathbf{v})$ or matrix $L$ corresponding to those elements we replace that row with the corresponding boundary operator or Fr\'echet derivative.

\begin{algorithm}[t]
	\caption{Overview of the adaptive  modular algorithm}\label{alg:overview}
	\begin{algorithmic}
		\State Initialize the geometry, boundary conditions, and algorithmic parameters
		\State Generate the initial guess $v$
		\While{res\_bvp $\geq$ tol\_bvp}
		\State Build $N(v)$ and $L(v)$
		\While{res\_new $\geq$ tol\_new}
		\State  $dv = -L(v)\backslash N(v)$	
		\State $v= v + dv$
		\State Compute res\_new
		\EndWhile
		\State Compute res\_bvp
		\If{res\_bvp $\geq$ tol\_bvp}
		\State Increase the number of Chebyshev and also Fourier points for cylindrical coordinates
		\State Sample the iterate solution to generate new initial guesses
		\EndIf
		\EndWhile
	\end{algorithmic}
\end{algorithm}

\section{Newton's method in the adaptive modular algorithm}
\label{newton}

We will treat the nonlinearity of $N$ with a Newton method.   The core of the algorithm is to first construct $N$ and $L$ as we have described above.  Then we solve the linear equation 
\begin{equation}
	L(\mathbf{v}) \, d\mathbf{v} = - N(\mathbf{v}).
\end{equation}
for  $d\mathbf{v}$ to build the update 
\begin{equation}
\mathbf{v}_{\mathtt{new}} = \mathbf{v} + d\mathbf{v}.
\end{equation}  
Here we solve for $d\mathbf{v}$ using Matlab's backslash operator, as even with our Kronecker products of dense matrices this is quite fast for almost all of the problems we have considered.
This iterative process continues until  the magnitude of the relative Newton residual $||d\mathbf{v}||/(||\mathbf{v}|| + \epsilon)$ is smaller than the prescribed tolerance, where we have set $\epsilon$ to $10^{-8}$ to avoid dividing by zero when the solution height is near 0.  
This process needs an initial guess to begin.  Generically we use either $\mathbf{v} = 0$ or $\mathbf{v} = 1$, though in some cases we have used more specific information about the solution when it is known beforehand.   If that initial guess is sufficiently close to the solution (and if some additional hypotheses hold) then the convergence is known to be quadratic.   For discussions of the functional analysis in Banach spaces and the associated convergence see Kantorovich and Akilov \cite{KantorovichAkilov1982}, as well as Ortega \cite{Ortega1968}, and later Gragg and Tapia  \cite{GraggTapia1974} gave optimal error bounds.  For a discussion of the convergence rate in finite dimensional problems see Nocedal and Wright \cite{NocedalWright2006}.  Here we do not check the additional hypotheses needed in these results, as that analysis is outside the scope of this work.

This inner loop is placed in a larger adaptive loop that checks the magnitude of the relative residual of the boundary value problem $||N(\mathbf{v})||/(||\mathbf{v}|| + \epsilon)$. If this residual is greater than the prescribed tolerance then this outer loop increases the number of points in the spectral method and resamples the iterate at this finer resolution to generate an initial guess to restart the inner Newton's method loop.   See Algorithm~\ref{alg:overview}. 

For this overview we have abstracted all of the details of the geometry and the details of the nonlinear and linear operators, as well as  the treatment of the  boundary conditions, given that we have defined these objects in a modular fashion in Section~\ref{operators}.


\section{Numerical results}
\label{results}

\begin{figure}[t!]
	\centering
	\scalebox{0.195}{\includegraphics{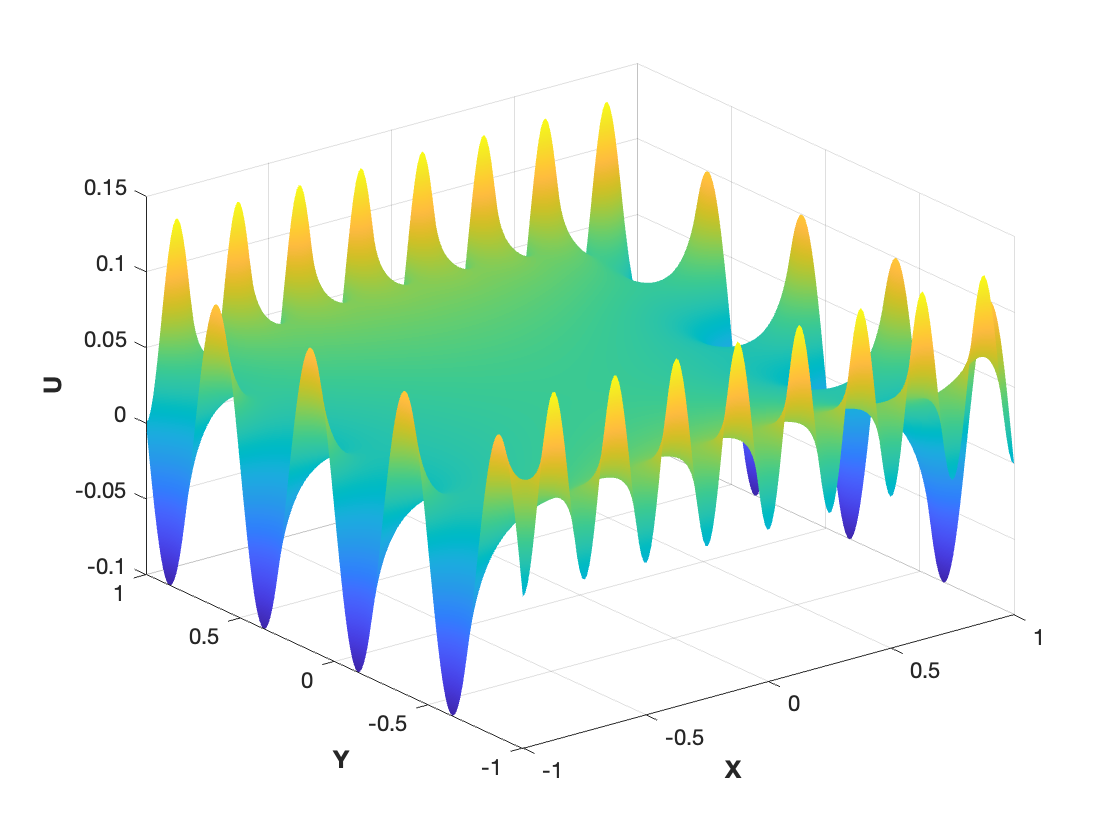}}
	\scalebox{0.195}{\includegraphics{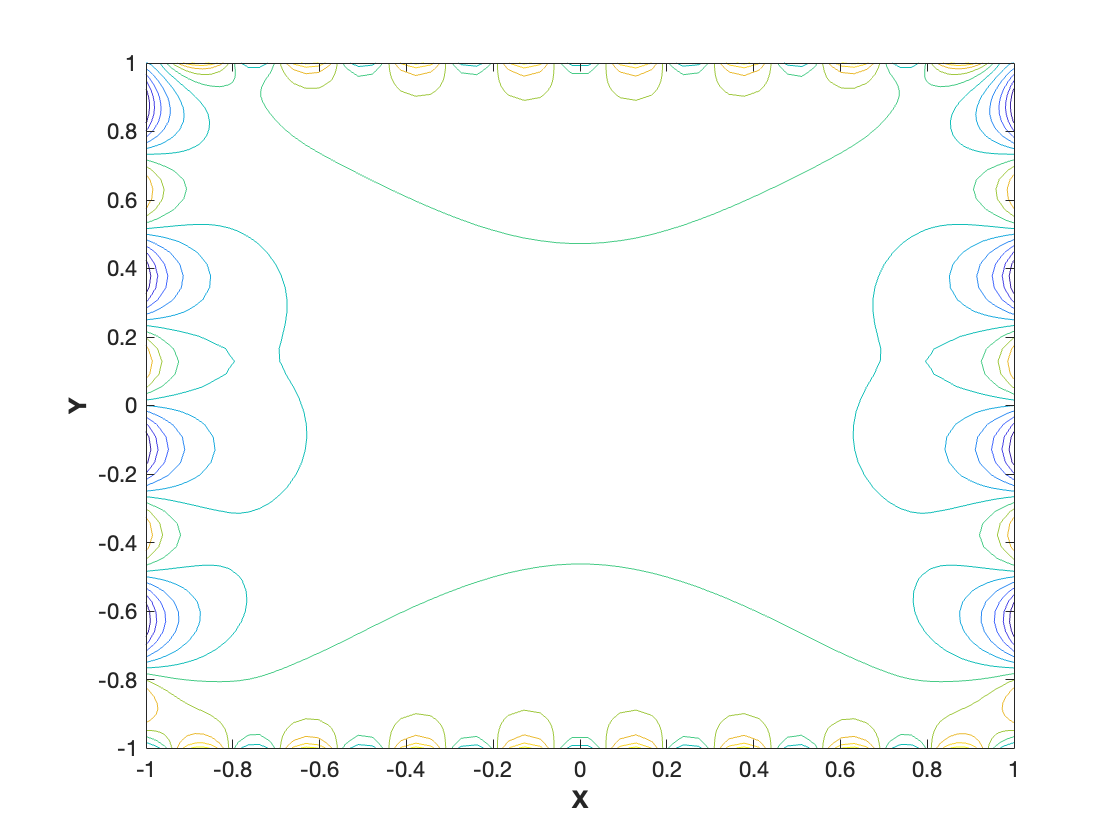}}
	\scalebox{0.195}{\includegraphics{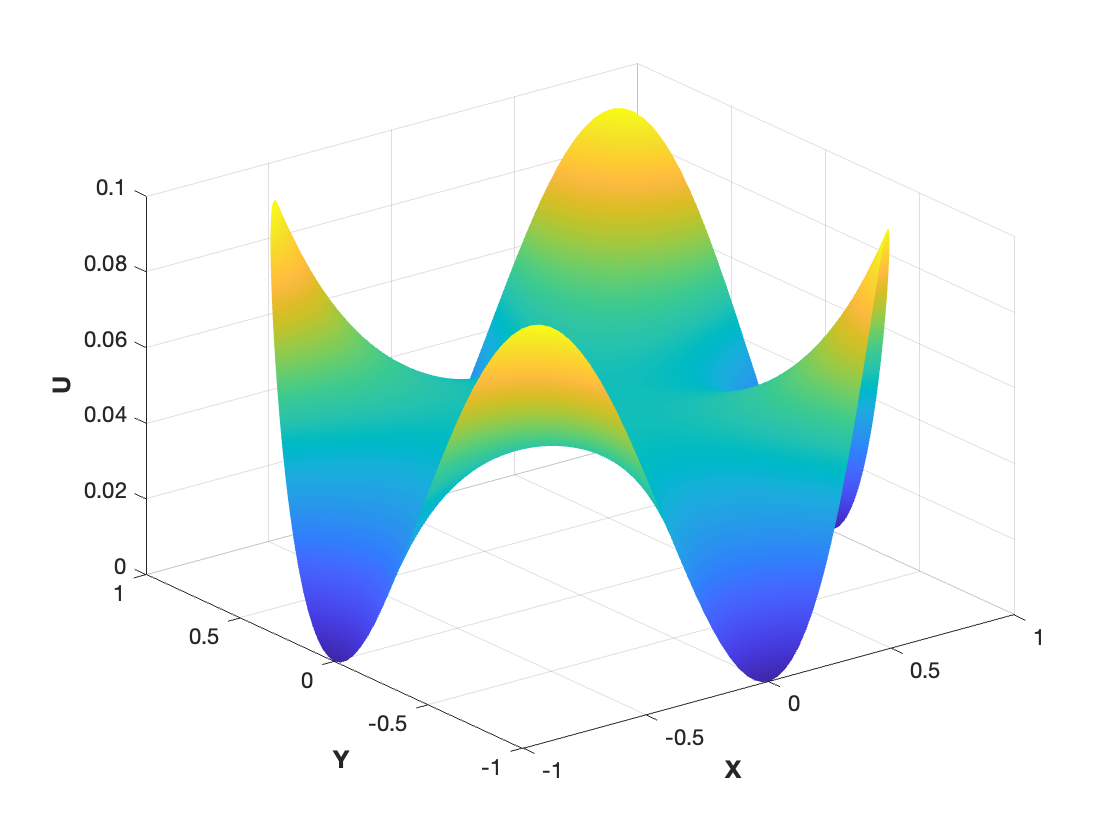}}
	\scalebox{0.195}{\includegraphics{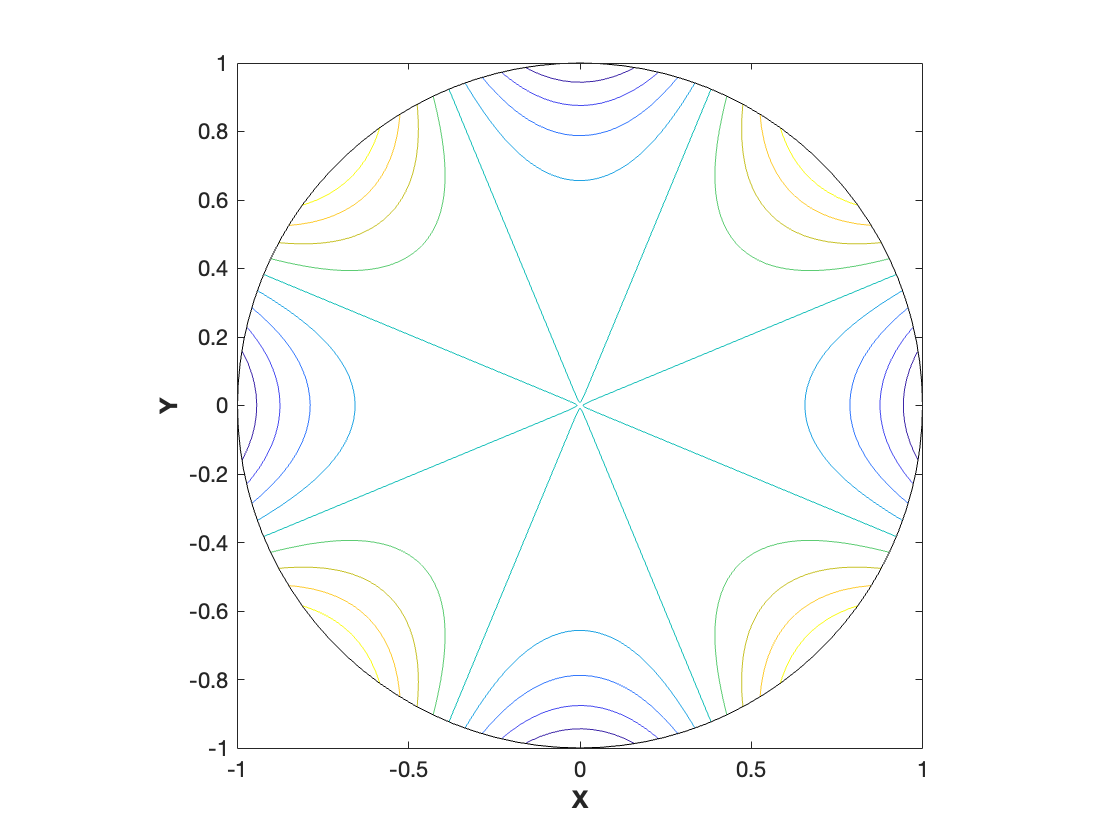}}
	\scalebox{0.195}{\includegraphics{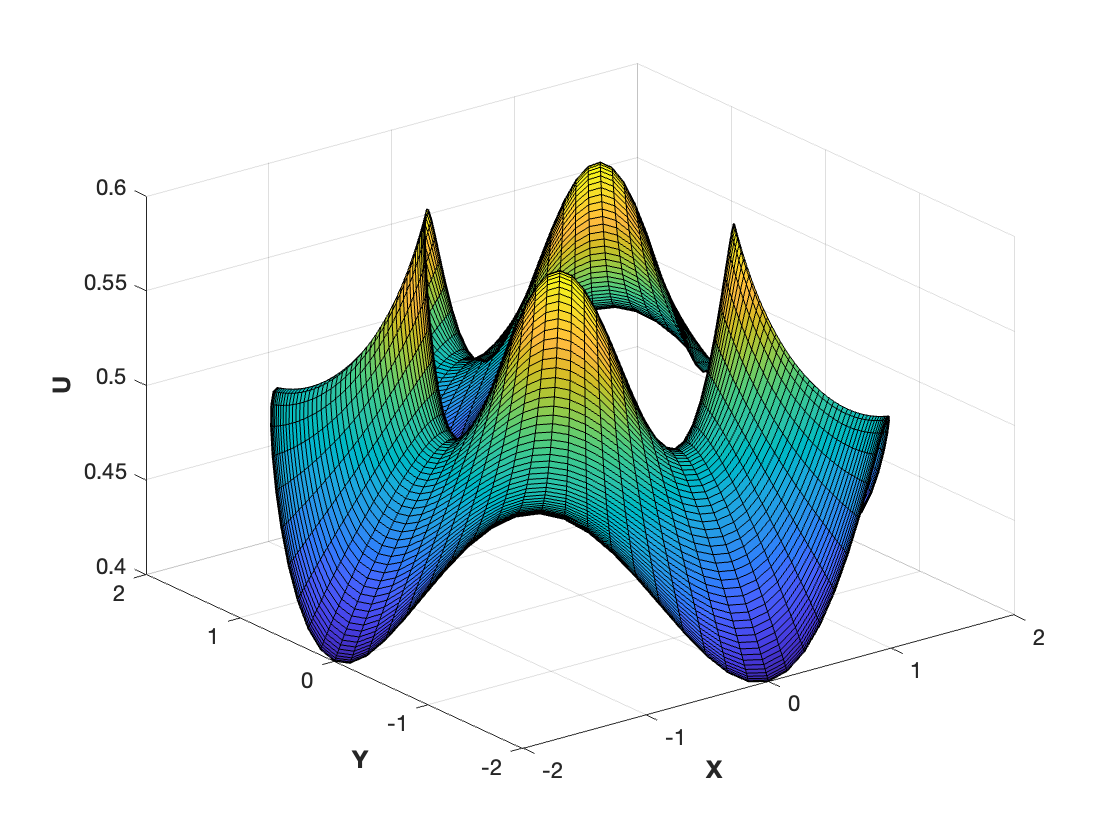}}
	\scalebox{0.195}{\includegraphics{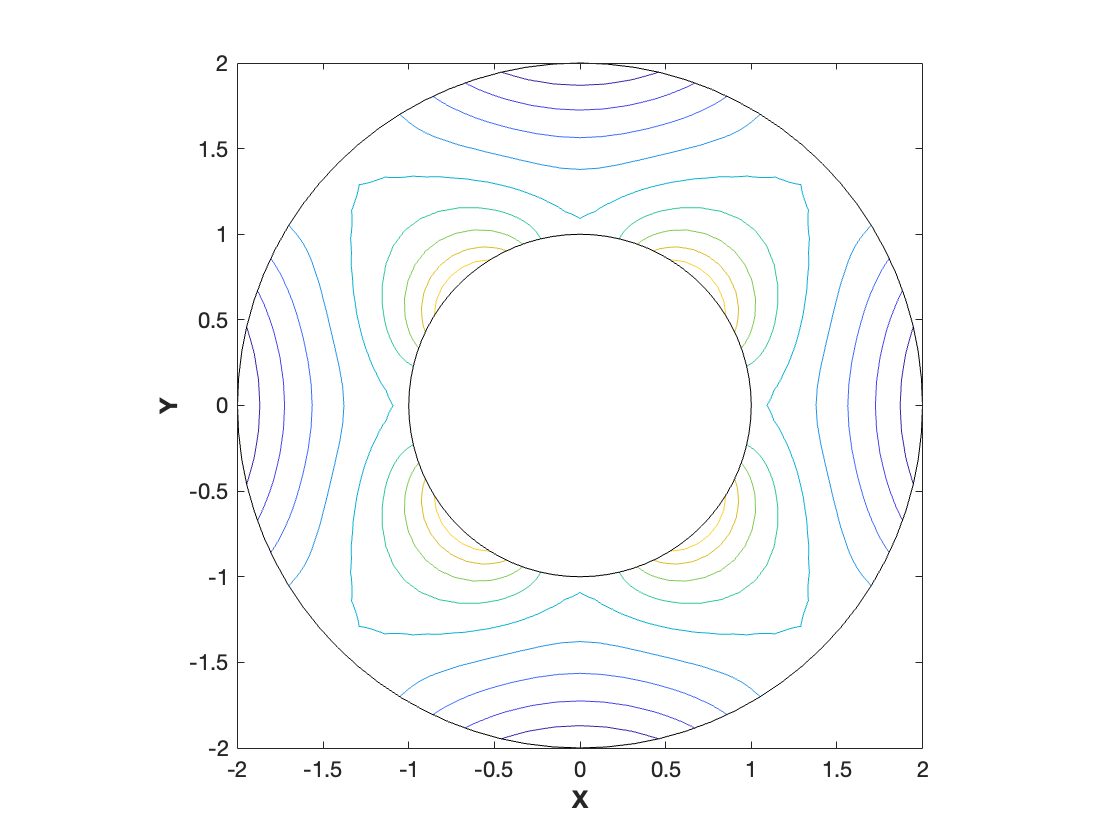}}
	\caption{The Plateau problem on three different domains.  The annular domain is not convex, so solutions do not exist for arbitrary boundary values.  }
	\label{fig:MinDir}
\end{figure}

\begin{figure}[h!]
	\centering
	\scalebox{0.195}{\includegraphics{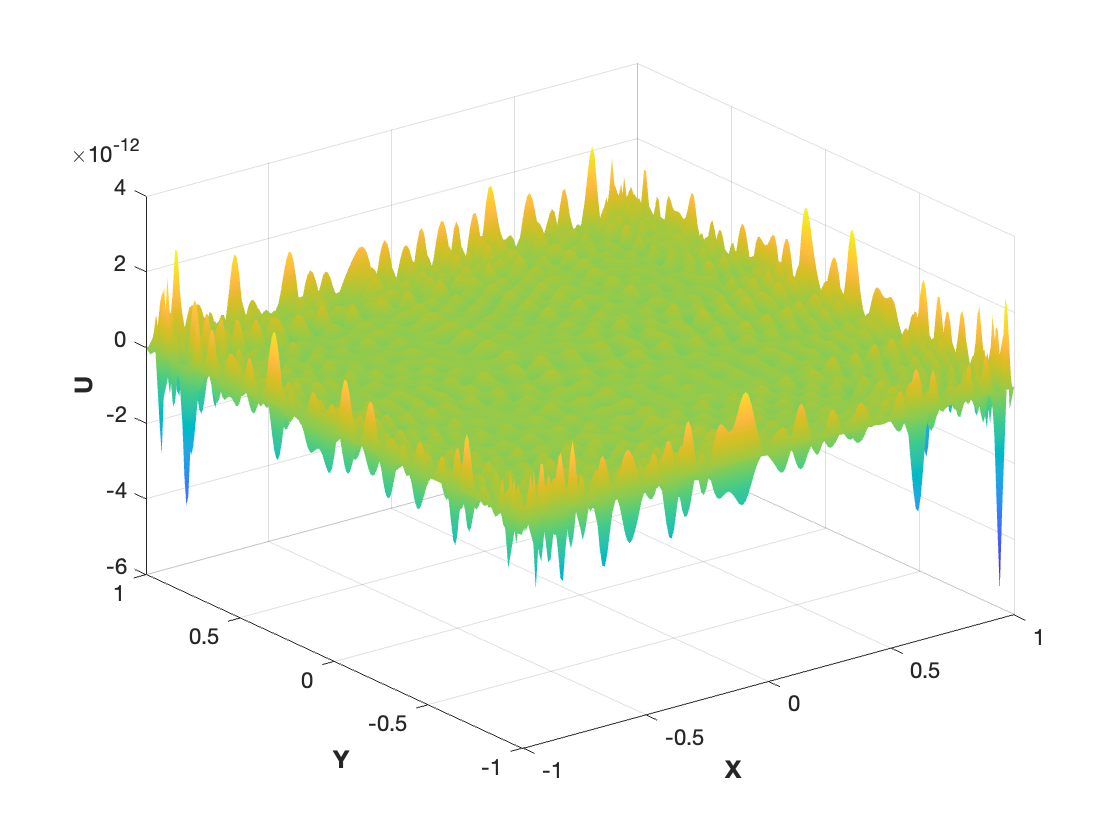}}
	\scalebox{0.195}{\includegraphics{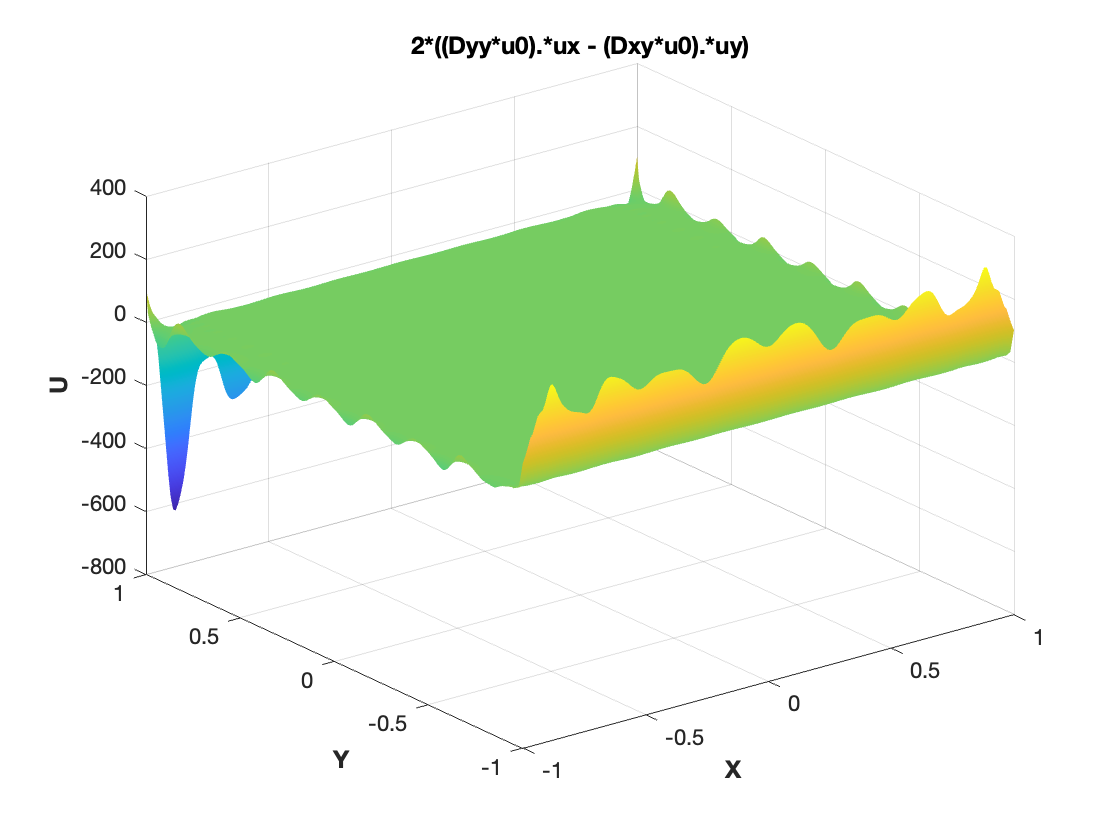}}
	\caption{On the left, the error in the solution of the Plateau problem for a minimal surface on the square with the number of Chebyshev points set to $n=55$.  Note that the vertical scale is 1e-12.  On the right we plot $2(u_{xx}u_y - u_{xy}u_x)$ which is the coefficient of the $D_x$ term in $F(u)$ for a convergent choice of coefficients given by $g(x,y)=  0.13\sin^2(4\pi x) + 0.1\sin(4\pi y)$.}
	\label{fig:MinDirRecErr}
\end{figure}

In the following experiments we choose a selection of Dirichlet boundary conditions and capillary boundary conditions as possible for minimal surfaces, CMC surfaces, and capillary surfaces over rectangles, disks, and annuli.  Our default tolerances for relative error are set to 1e-14 for the norm of the Newton step on the rectangle, and 1e-13 for the polar forms of the problems.  The default tolerance for the relative error of the boundary value problem as measured by $||N(u)||$ is 1e-10 for all of the problems we consider.  The default number of Chebyshev points we use for square domains is $n = 55$.  We arrived at this default number for most of our problems by observing plots of $N(u)$ as the number of grid points increased.  These error plots stabilized with maximum magnitudes of the error on the order of 1e-12 at $n=55$ and the errors were somewhat erratic for coarser resolutions.  We show this in Figure~\ref{fig:MinDirRecErr} for the Plateau problem for minimal surfaces shown at the top of Figure~\ref{fig:MinDir}. We will discuss that example shortly.  We prefer to use these default resolutions rather than appealing directly to the adaptive portions of the algorithm because the resulting code runs much faster.  However, the adaptive portions of the algorithm are particularly useful if the domain is modified from our standard configuration.   For the Chebyshev-Fourier methods, we used defaults of $n=50$ Chebyshev points and $m=80$ Fourier points, and we will discuss our reasons for that ratio below.  For the disk domains we do have a double covering so that resolution is essentially halved, and we keep the higher resolution for the annular problems.

We begin with the figures earlier in the paper.  Figure~\ref{fig:MinDirAnnCat} shows a minimal surface with Dirichlet boundary data with a domain of an annulus.  The height of the surface at the outer radius is prescribed to be 0 and at the inner radius it is prescribed to be a positive value.  This value was chosen to get reasonably close to the gradient blow-up.   The initial guess for the height was 0.  The plotting commands in Chebfun have not yet been extended to annular domains, so we merely used the numerical grid points to plot the solution.  The actual interpolant has a much smoother character.  We keep the mesh grid as a visual cue for these annular problems.   In fact, we chose more angular points simply because the resulting graphics were smoother, and the algorithm converged with the same tolerances when $m=50$ Fourier points were used.  In contrast, with Figure~\ref{fig:CapCorner} we have the plotting commands in Chebfun, and we are able to use barycentric interpolation to more accurately represent the polynomial interpolant on that domain.  This smoother plotting will be used whenever possible.  The numerical experiment was to approximate the capillary surface with contact angle close to the optimal value, and we got to 0.035 above the optimally extreme value of $\pi/4$.  We used $n = 50$ for the number of Chebyshev points and the default tolerances.  We used both initial guesses of the height identically 1 and the height defined by 
$$
v_0 = 2\sqrt{2} + \frac{1}{2\sqrt{2}} - \sqrt{8 - x^2 - y^2}
$$
and we found no meaningful difference in the convergence times.  This particular guess is based on estimates found in \cite{ecs}.  The loops converged in 7.774170 seconds for the constant guess and 7.771485 seconds for the guess that uses some prior knowledge of the solution.  These numbers come from using a Macbook Pro from 2023 with an M2 chip, and we use that machine for all measurements of elapsed time in this paper.

Figure~\ref{fig:MinDir} shows examples of the minimal surfaces and the corresponding Plateau problem with different domains.  The top figure with the square domain shows an example with boundary conditions given by the function
$$
g(x,y) =  0.1\sin^2(4\pi x) + 0.1\sin(4\pi y)
$$
to make an interesting figure.  The coefficients were chosen so that the algorithm converges.  If the amplitudes are increased by much then the matrix $L(v)$ becomes ill conditioned.  For example, if we set $g(x,y) =  0.13\sin^2(4\pi x) + 0.1\sin(4\pi y)$ the algorithm does converge, but much larger coefficients do not give convergence.  
In exploring the limits of the variation of $g$ we also coded a method of continuation.  The approach we took was to multiply $g$ by a parameter $\lambda$ and then start out with an initial guess for the solution $u\equiv 0$ when $\lambda = 0$.  Then we used the converged solution for that $\lambda$ value as an initial guess for a slightly increased value for $\lambda$.  This was done for a sequence of values up to $\lambda = 1$.  Our conclusion from this experiment is that the limits on the variation of $g$ are mathematical, and not based on poor initial guesses.  We illustrate what is happening in the right side of Figure~\ref{fig:MinDirRecErr}, which shows the coefficient $2(u_{yy}u_x - u_{xy}u_y)$ for the $D_x$ term in $F(u)$ with a magnitude on the order of 600 along the boundary.  The $D_y$ coefficient has similar magnitudes along the boundary as well.  In contrast, the amplitude of 0.1 leads to the coefficient of $D_x$ having a magnitude of 200 along the boundary and while the coefficient of $D_y$ some regions with a magnitude of 600, these regions are much smaller and isolated at the corners.  This is what guided our limitations on the coefficients in $g$.  For the example in our figure, we used $n = 86$ Chebyshev points in an effort to smooth the contour plot.  The middle figure with the disk domain has boundary conditions given by
$$
g(\theta) = 0.1\sin^2(2\theta)
$$
and default settings were used.  The bottom figure has an annular domain with inner radius $a = 1$ and outer radius $b = 2$.  The boundary conditions are
$$
g(\theta) = 0.5 + 0.1\sin^2(2\theta)
$$ 
at $r = a$ and 
$$
g(\theta) = 0.5 - 0.1\cos^2(2\theta)
$$
at $r = b$ and default settings were used.

\begin{figure}[h!]
	\centering
	\scalebox{0.195}{\includegraphics{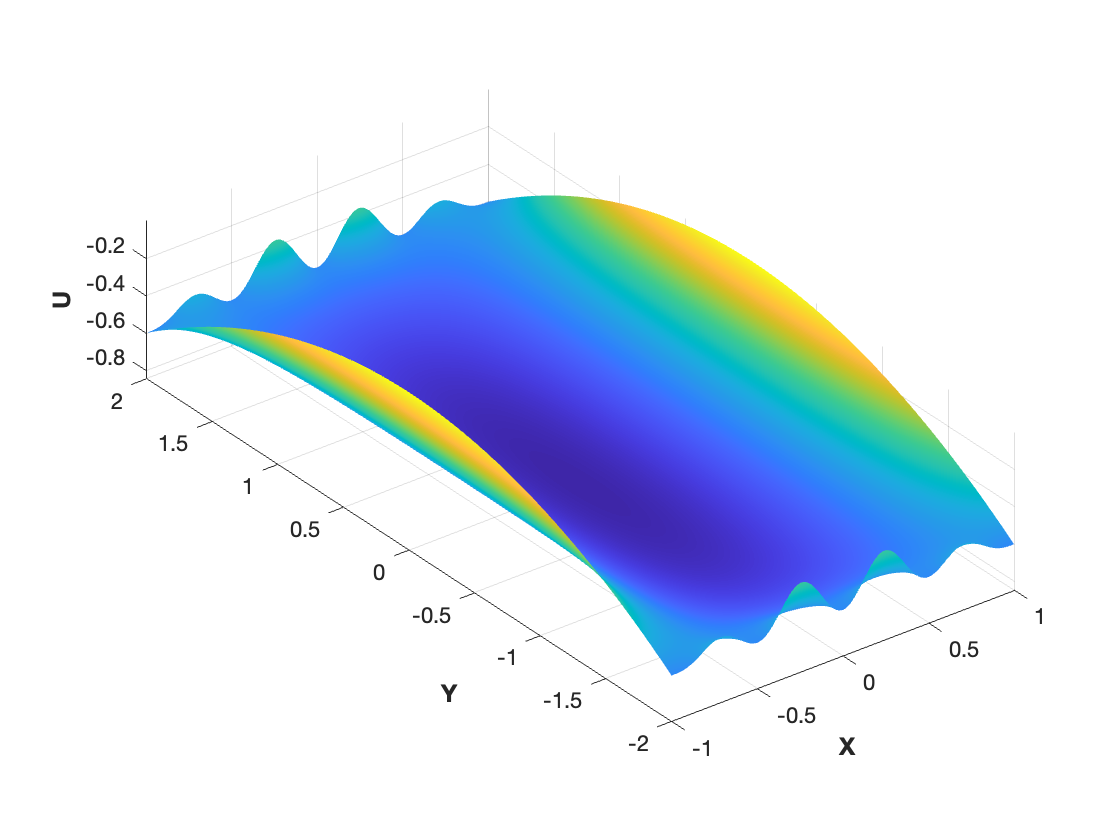}}
	\scalebox{0.17}{\includegraphics{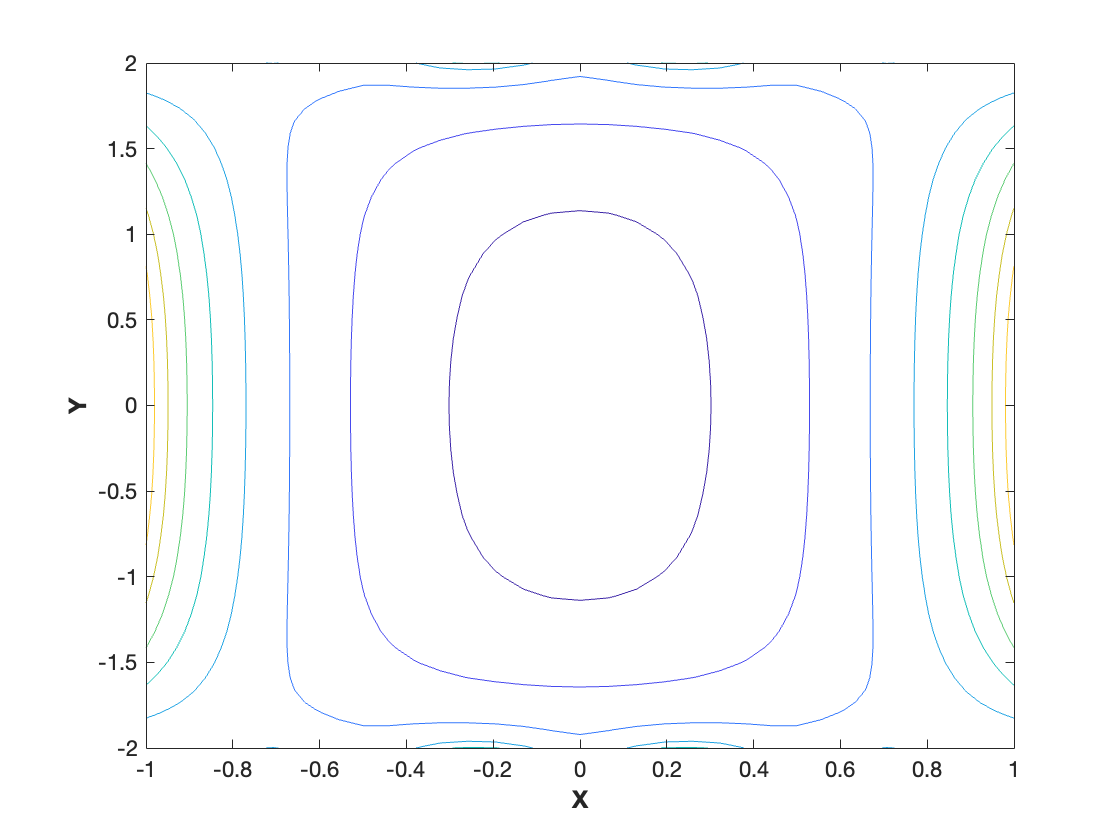}}
	\scalebox{0.195}{\includegraphics{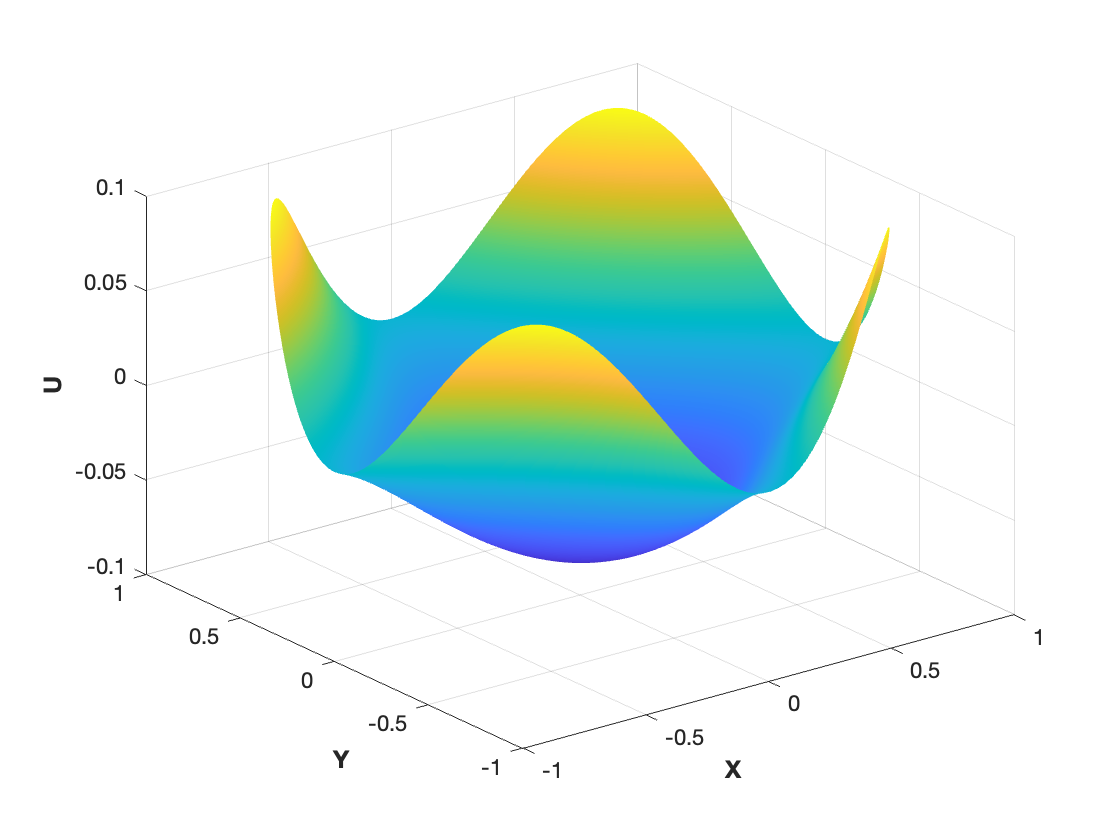}}
	\scalebox{0.195}{\includegraphics{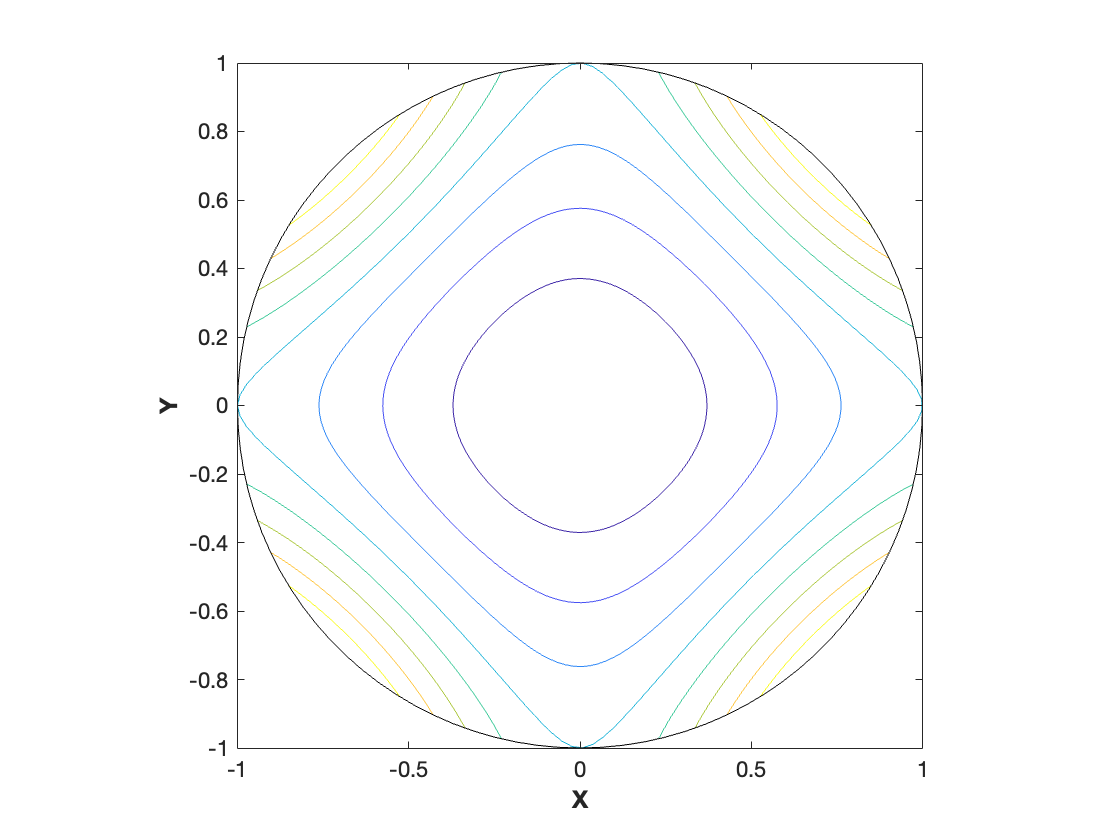}}
	\scalebox{0.195}{\includegraphics{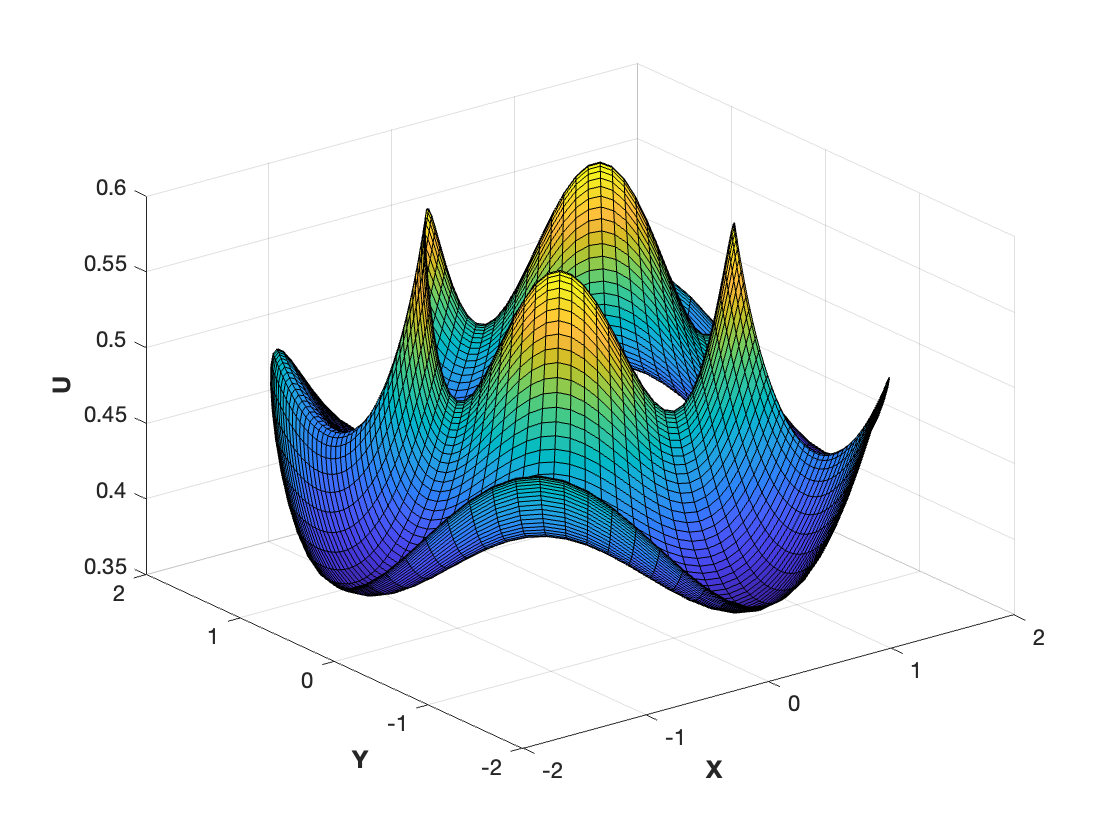}}
	\scalebox{0.35}{\includegraphics{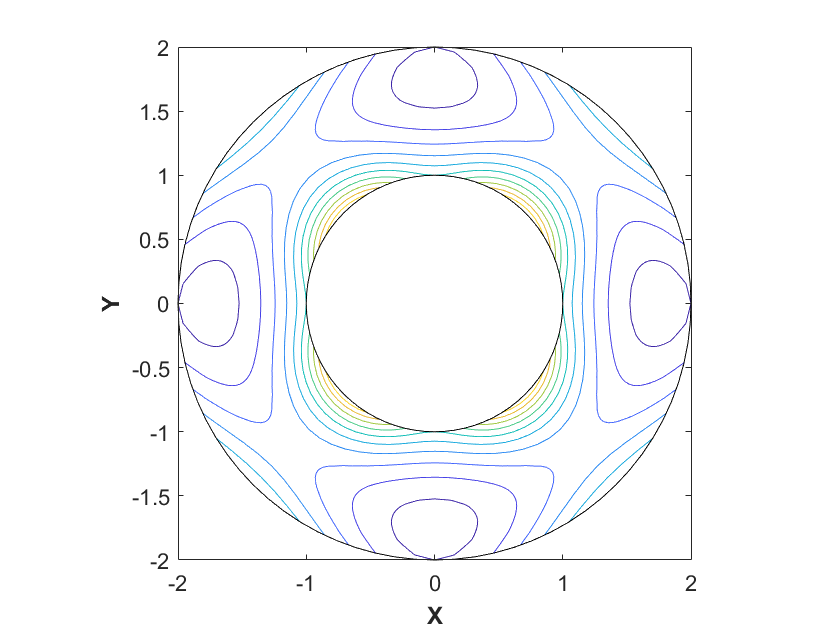}}
	\caption{The Plateau problem for CMC surfaces on three different domains.   The annular domain is not convex, so solutions do not exist for arbitrary boundary values. }
	\label{fig:CMCDir}
\end{figure}

\begin{figure}[h!]
	\centering
	\scalebox{0.195}{\includegraphics{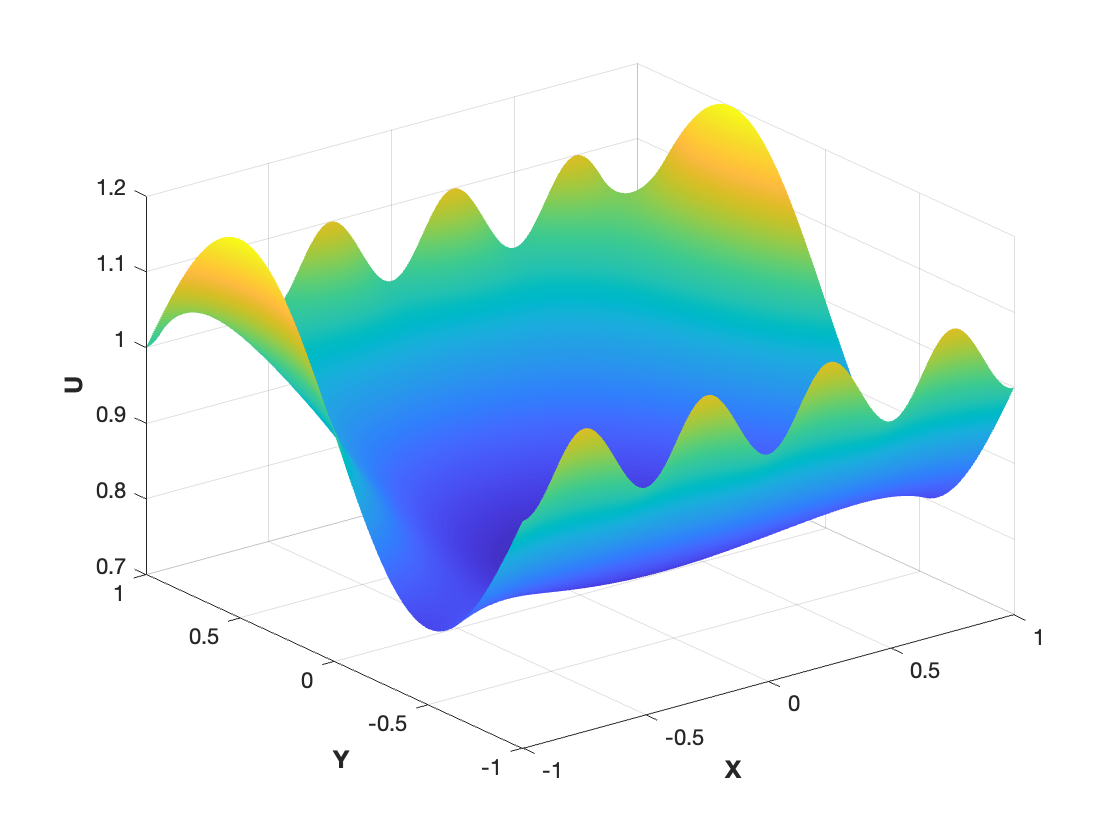}}
	\scalebox{0.195}{\includegraphics{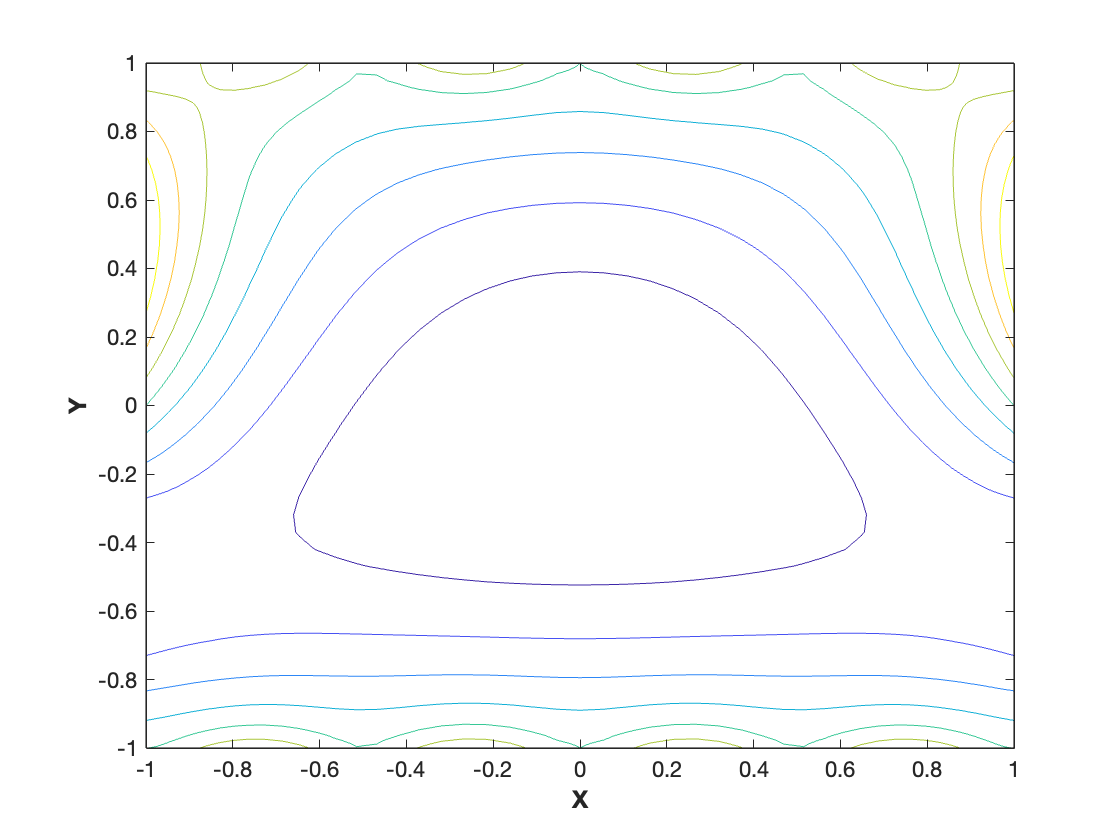}}
	\scalebox{0.195}{\includegraphics{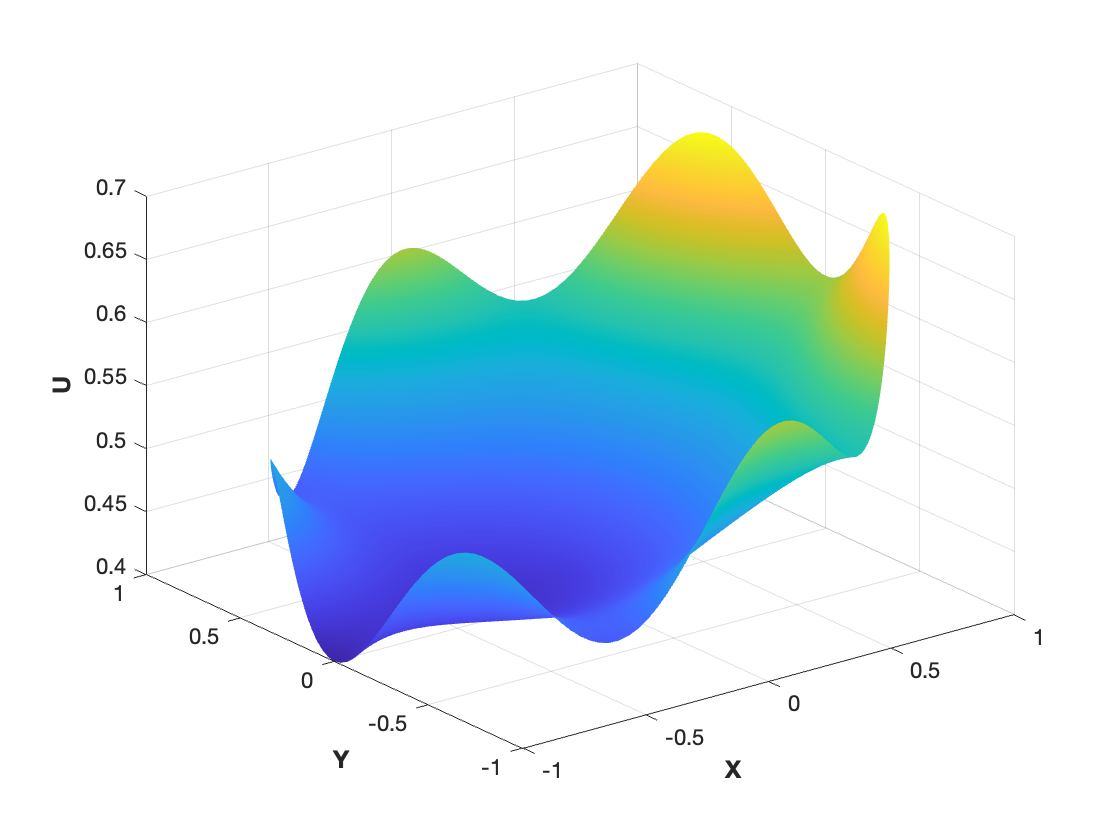}}
	\scalebox{0.195}{\includegraphics{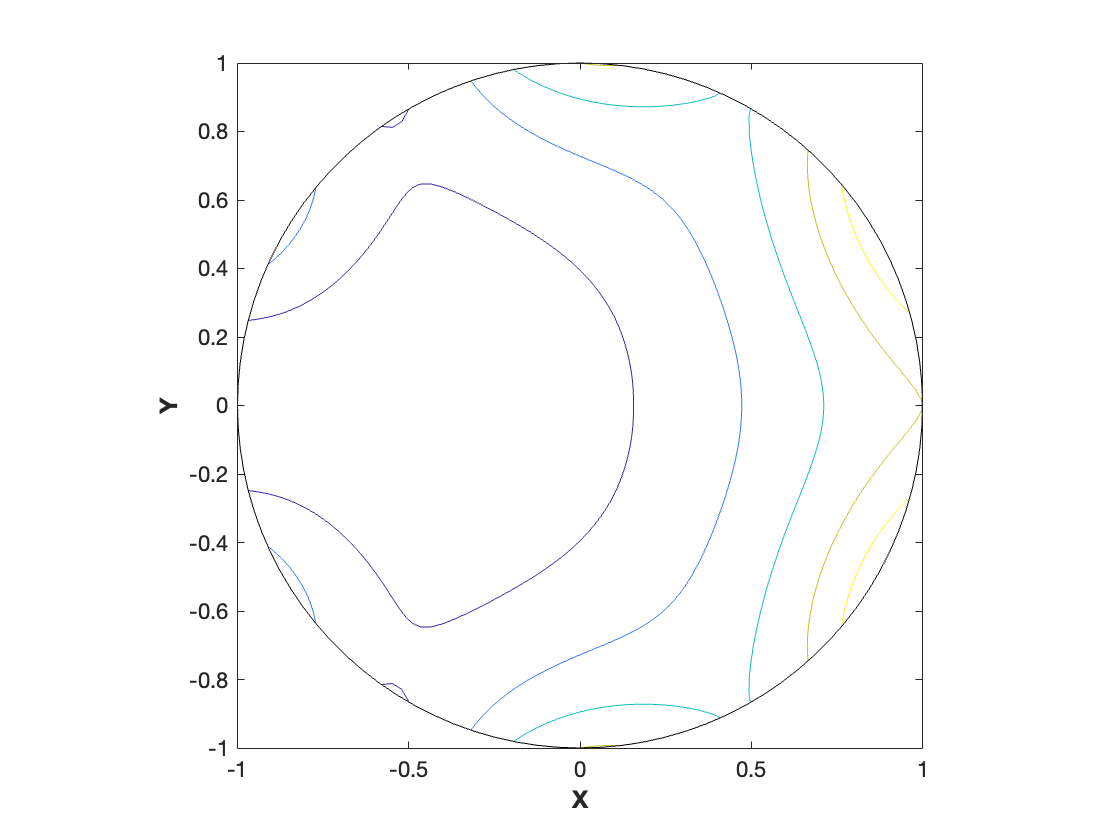}}
	\scalebox{0.195}{\includegraphics{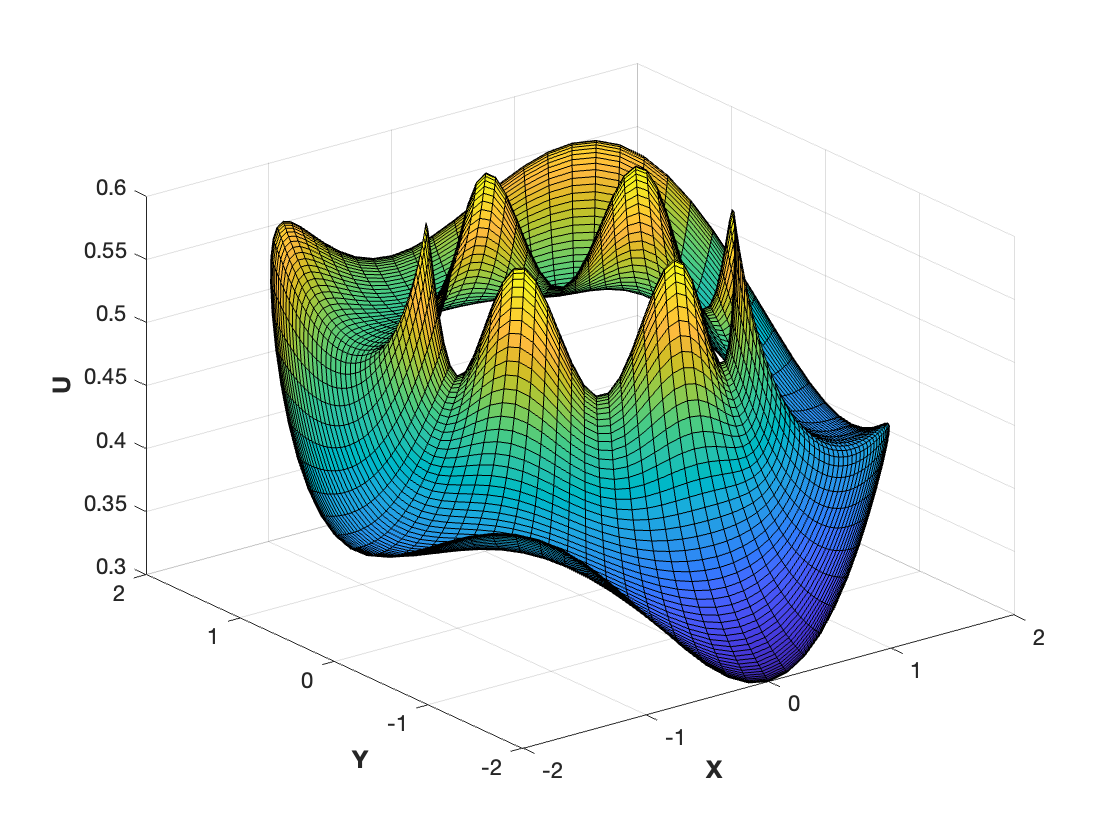}}
	\scalebox{0.195}{\includegraphics{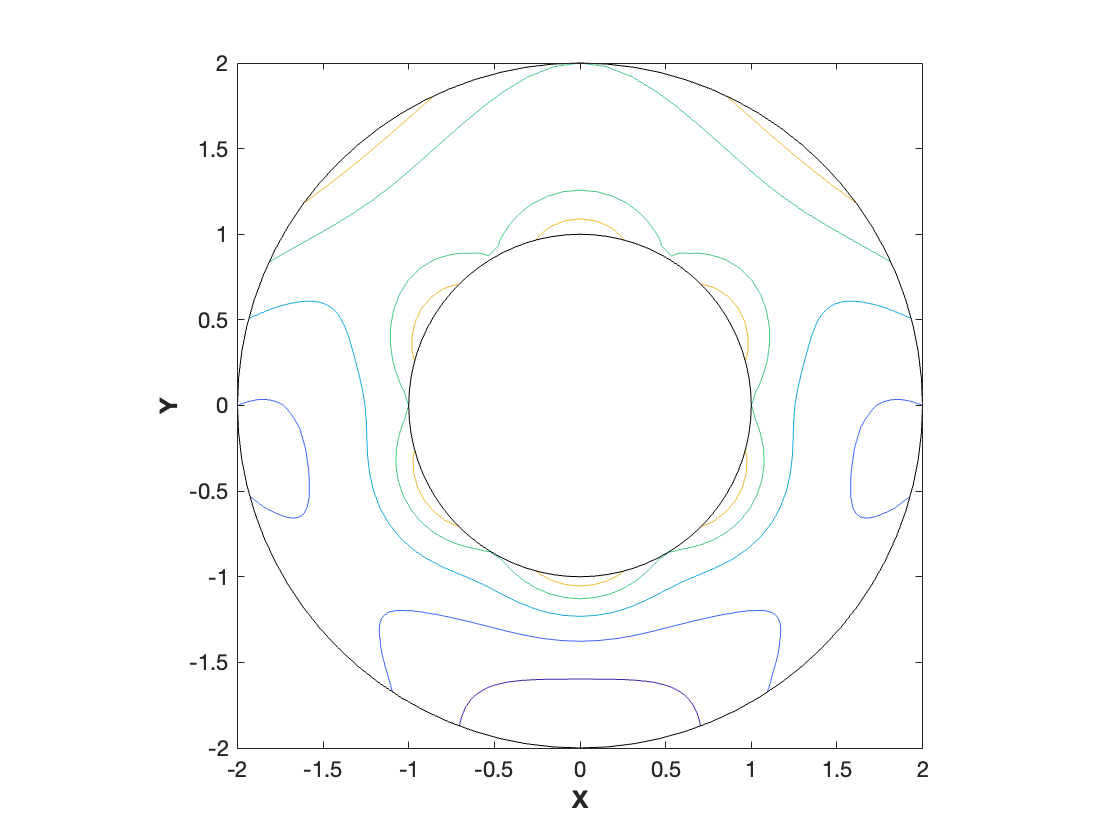}}
	\caption{The Plateau problem for capillary surfaces on three different domains.   }
	\label{fig:CapDir}
\end{figure}

Figure~\ref{fig:CMCDir} shows examples of the Plateau problem for CMC surfaces over different domains.  The top figure with a rectangular domain of $[-1,1]\times[-2,2]$ has boundary conditions given by the function
$$
g(x,y) =  0.25\cos(\pi x/2)\sin^2(2\pi x) - 0.15y^2
$$
with $2H = 1$.    The middle figure with a disk domain has $2H = 1/2$ and boundary data 
$$
g(\theta) = 0.1\sin^2(2\theta).
$$
like the minimal surface in Figure~\ref{fig:MinDir}.   For the annular domain in the bottom figure we also repeated the experiment from the minimal surface in the previous figure, with $2H = 1/2$ again.

Figure~\ref{fig:CapDir} shows examples of the capillary surfaces and the corresponding Plateau problem with different domains.  The top figure with a rectangular domain of $[-1,1]\times[-1,1]$ has boundary conditions given by the function
$$
g(x,y) =  0.1\sin^2(2\pi x) + 0.2\sin(\pi y) + 1;
$$
with $\kappa = 1$.  To achieve the tolerance of 1e-14 the adaptive algorithm used $n = 59$ Chebyshev points.  The middle figure with a disk domain has boundary data 
$$
g(\theta) = 0.1\sin^2(3 \theta) + 0.1\cos(\theta) + 1/2.
$$
For the annular domain in the bottom figure the boundary data was
$$
g_a(\theta) = 0.5 + 0.1\sin^2(3\theta),
$$
and
$$
g_b(\theta) =  0.5 - 0.1\cos^2(2\theta) + 0.1\sin(\theta)
$$
on the inner and outer radii of $a = 1$ and $b = 2$ respectively.

\begin{figure}[h!]
	\centering
	\scalebox{0.195}{\includegraphics{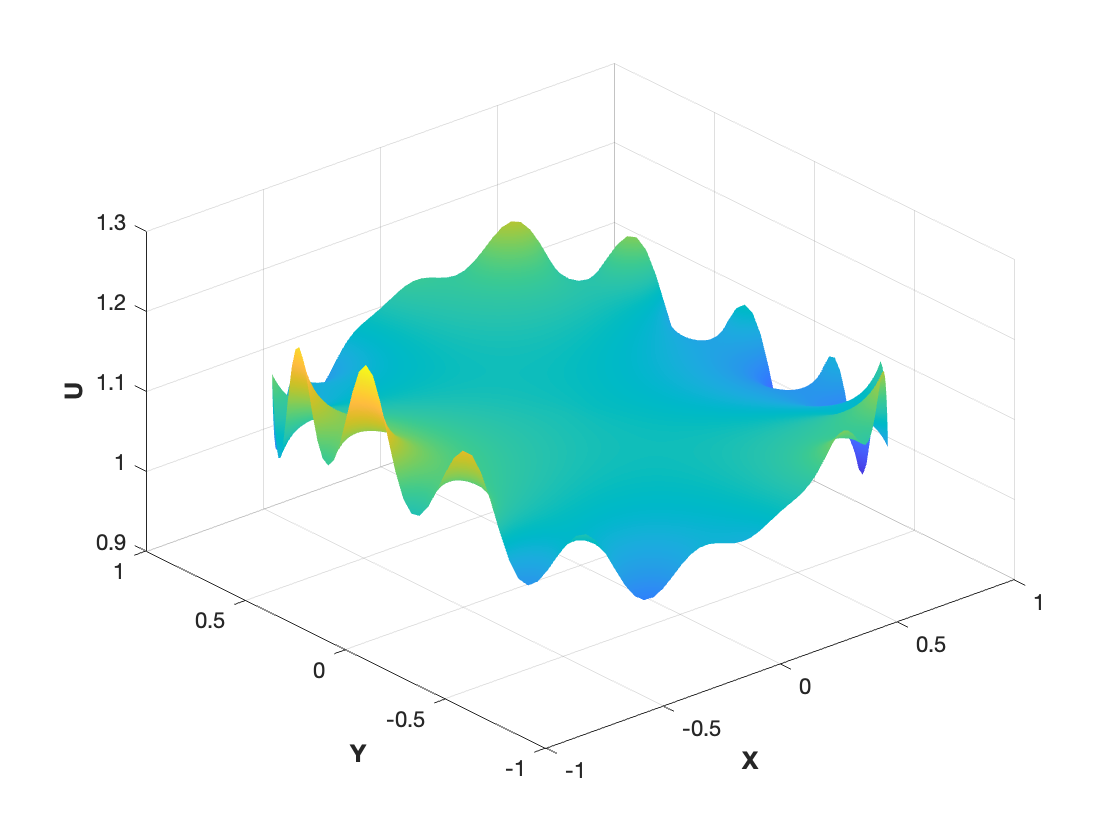}}
	\scalebox{0.18}{\includegraphics{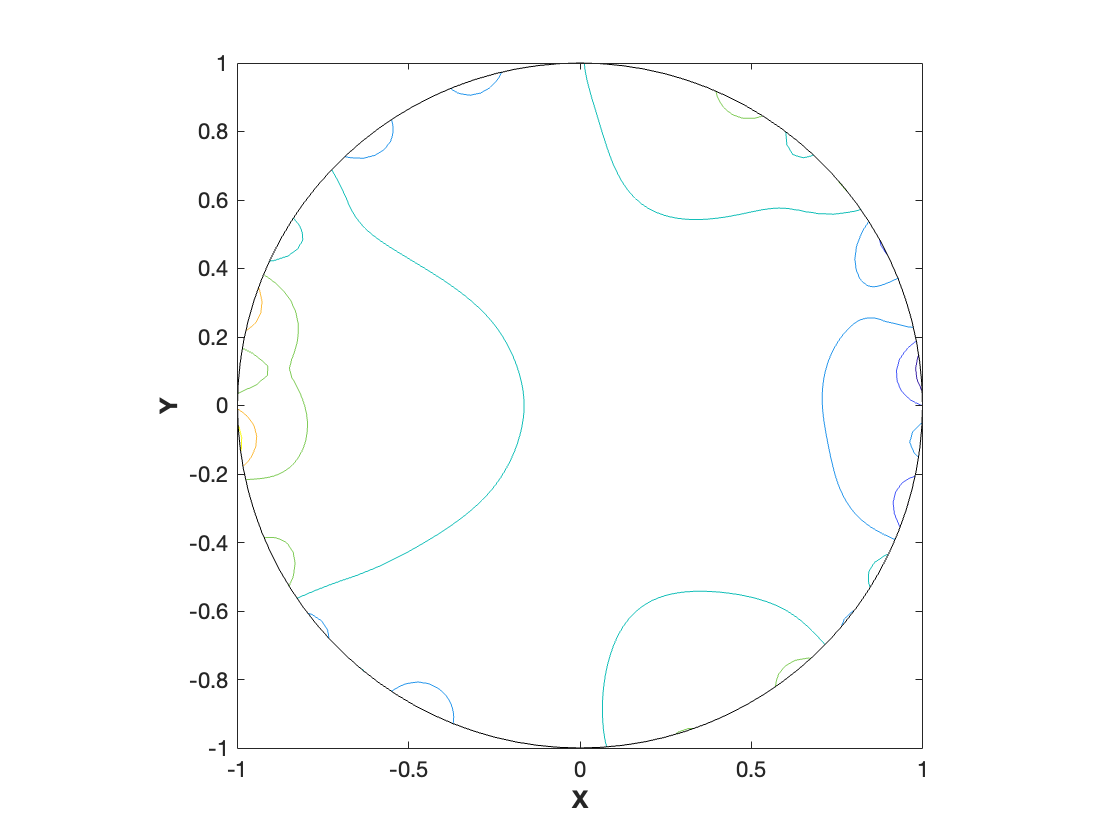}}
	\scalebox{0.195}{\includegraphics{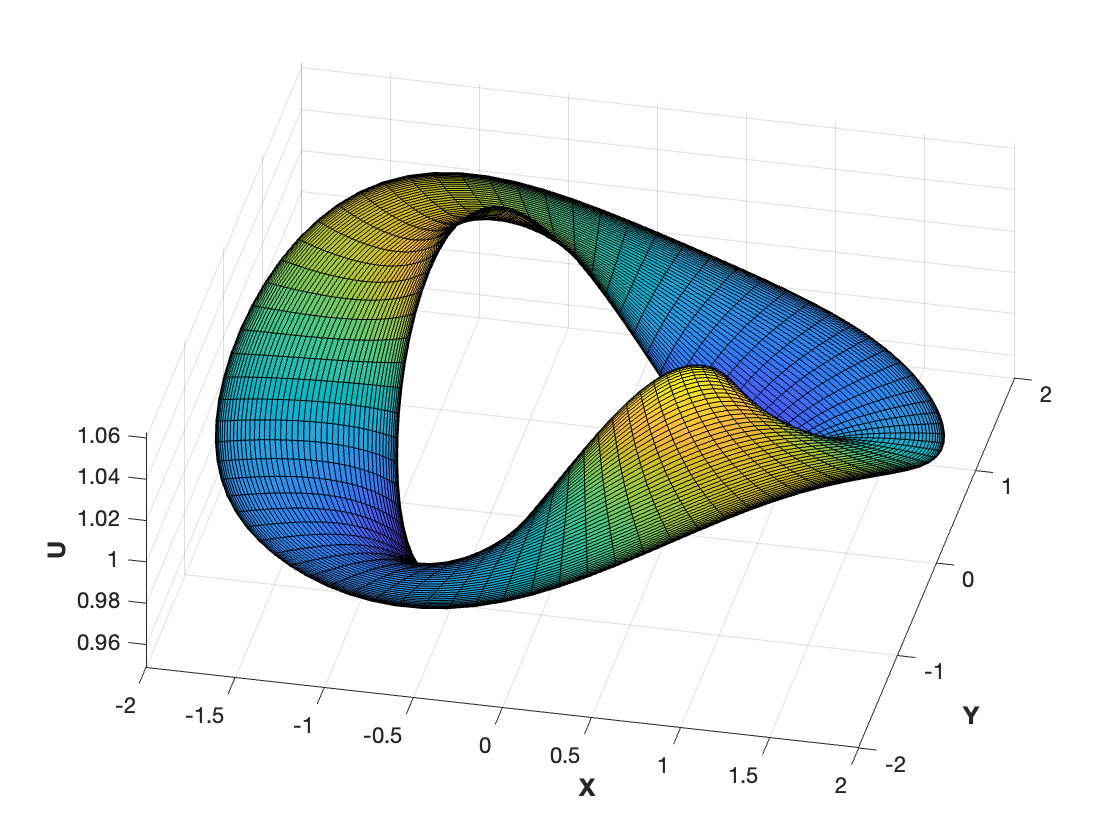}}
	\scalebox{0.18}{\includegraphics{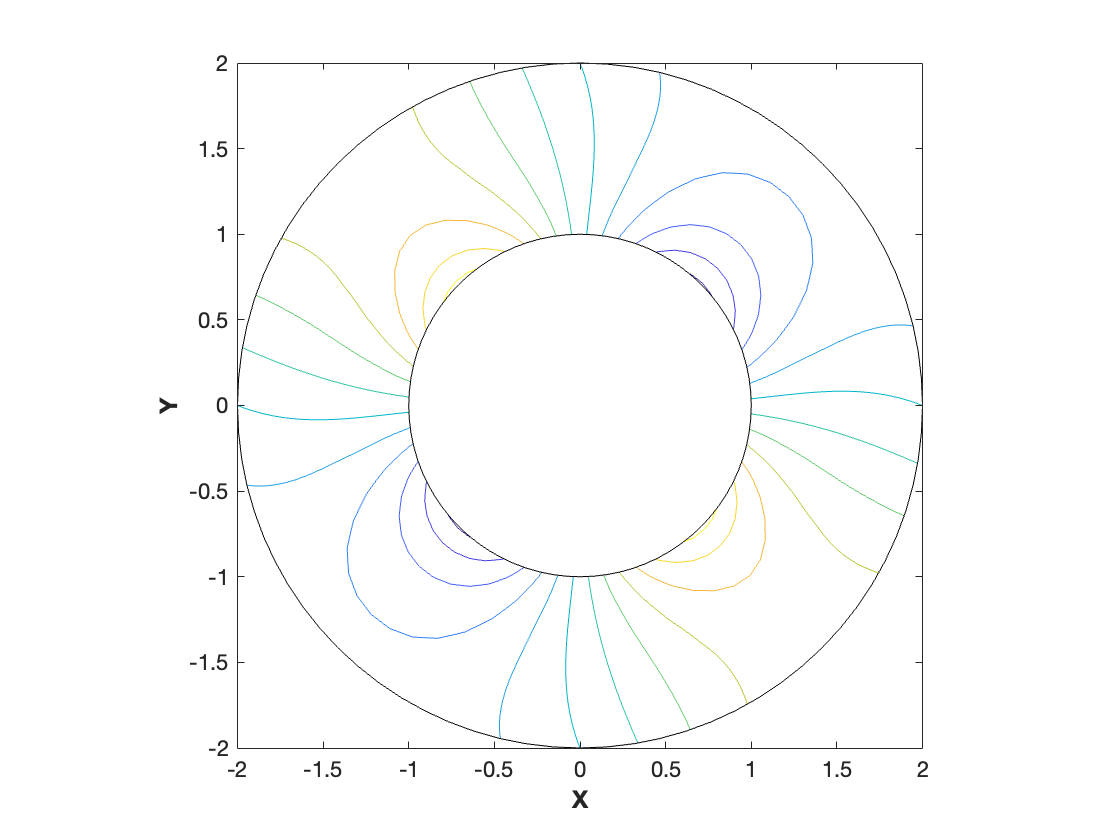}}
	\caption{Minimal surfaces with capillary data on two different domains.}
	\label{fig:MinNeu}
\end{figure}

For both the minimal surfaces with capillary boundary data, and for the CMC surfaces with the corresponding capillary problems without gravity our algorithms did not converge for of the rectangular domains we considered.  
The mean curvature is determined by the domain and the boundary conditions here, as we noted in \eqref{eqn:2H}.  The solutions to these equations are only unique up to vertical translation, so one can prescribe the height at one point to give a unique solution.  We implemented this for the rectangular domain and our algorithm did not converge.  We also used Robin boundary conditions at the corner points to average the capillary boundary conditions with prescribed heights there and that did not converge.  Finally we tried additionally prescribing the average of the height over the entire boundary and this also did not converge.  We attribute our codes not working to the lack of stability of these surfaces as in Theorem~\ref{CMCunstable}.  

In contrast, we were able to get our algorithms to converge for disk and annular domains.  Generically, we specified a height of 1 for the $\theta = 0$ angle on the boundary of the disk or on the outer boundary of the annulus.   This choice of the height at one point leads us to choose an initial guess with uniform height of 1.

Figure~\ref{fig:MinNeu} shows minimal surfaces with capillary boundary data.  We require  $\int_{\partial\Omega}\cos(\gamma(s))\, ds = 0$.  On the top we have a disk domain with
$$
\gamma(\theta) = \frac{\pi}{2}  +  \sin(16\theta)\cos(\theta) + 0.2\cos(3\theta) + 0.05\cos(\theta).
$$
On the bottom we have an annular domain with $a = 1$ and $b = 2$.  We used
$$
\gamma_a = \frac{\pi}{2} + 0.1\sin(2\theta)
$$
and
$$
\gamma_b = \frac{\pi}{2} + 0.01\cos(4\theta)
$$
at $a$ and $b$ respectively.  We note that our algorithm is somewhat sensitive to the boundary data, so the coefficients we have were chosen carefully.  Also, we used an initial guess of uniform height 0 as that performed better than one with uniform height 1.

\begin{figure}[h!]
	\centering
	\scalebox{0.195}{\includegraphics{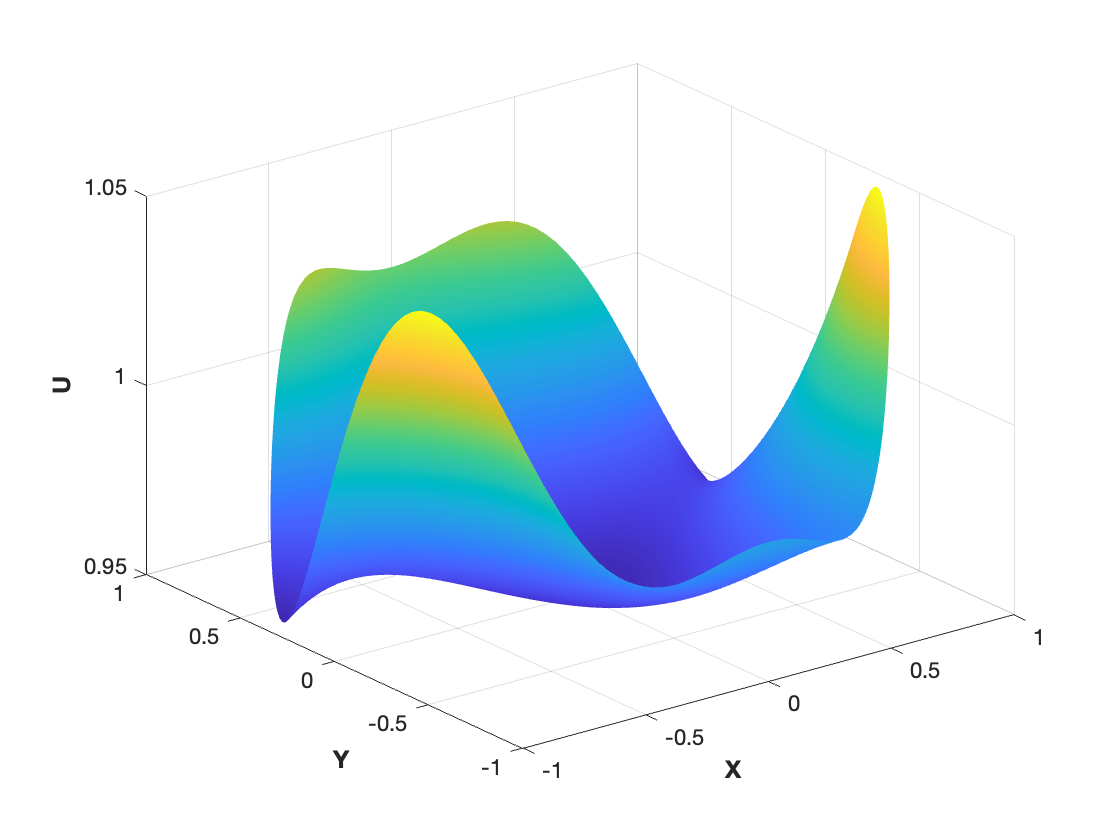}}
	\scalebox{0.18}{\includegraphics{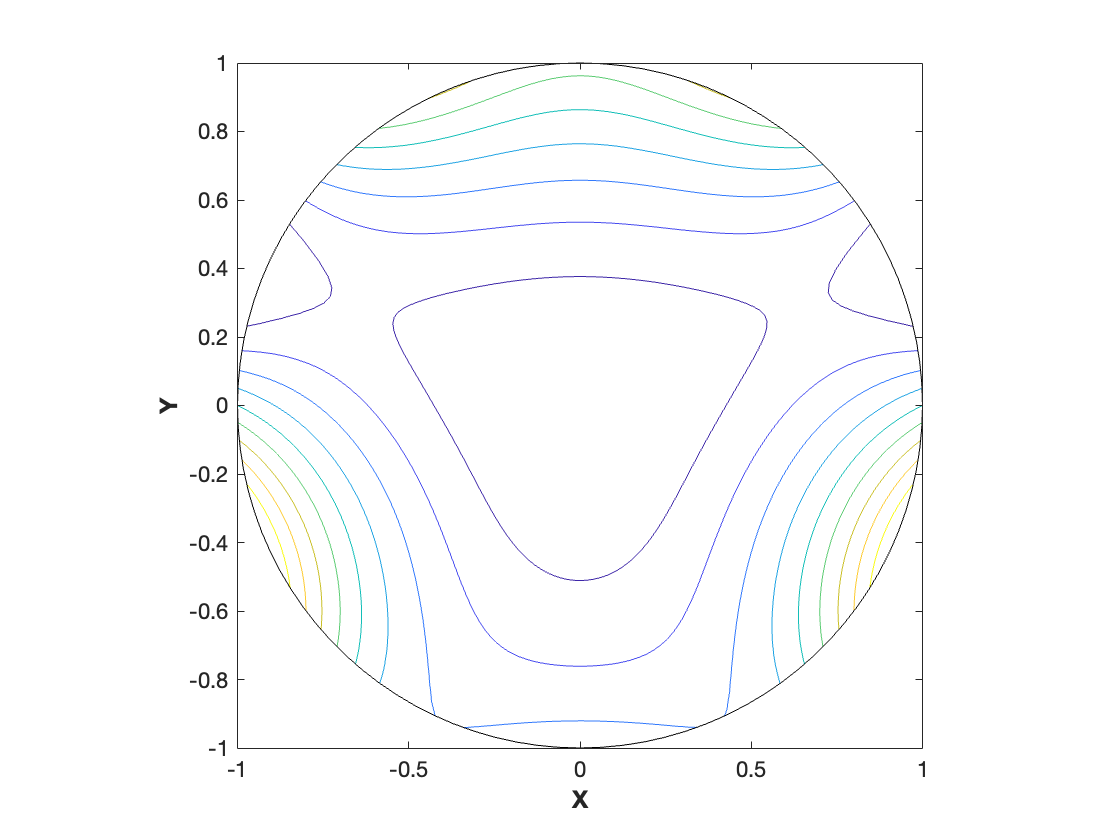}}
	\scalebox{0.195}{\includegraphics{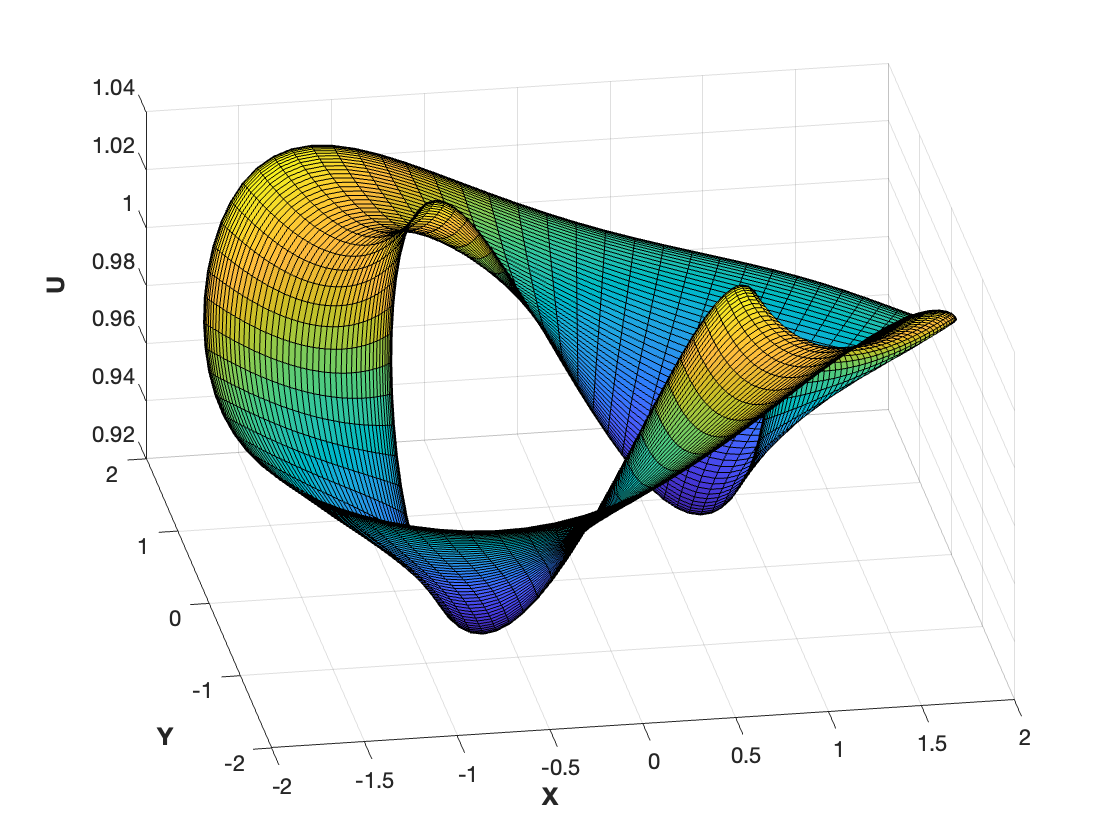}}
	\scalebox{0.18}{\includegraphics{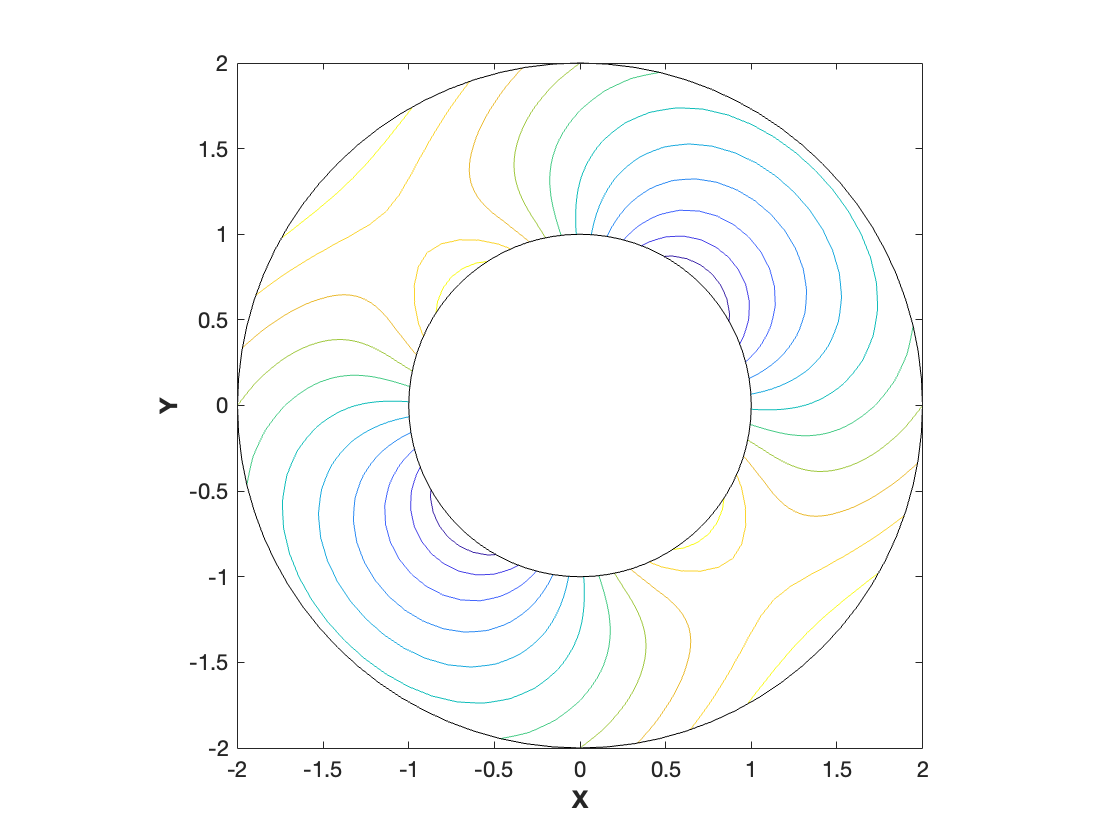}}
	\caption{Zero gravity capillary surfaces, or CMC surfaces with capillary data on two different domains.  }
	\label{fig:CMCNeu}
\end{figure}

Figure~\ref{fig:CMCNeu} shows CMC surfaces with capillary boundary data, which are also capillary surfaces without gravity.  On the top we have a disk domain with
$$
\gamma(\theta) = \frac{\pi}{2}  - 0.1 +  0.2\sin(4\theta)\cos(\theta),
$$
which leads to $2H =  0.1987$.
On the bottom we have an annular domain with $a = 1$ and $b = 2$.  We used
$$
\gamma_a = \frac{\pi}{2} + 0.1\sin(2\theta)
$$
and
$$
\gamma_b = \frac{\pi}{2} - 0.05 + 0.01\cos(4\theta)
$$
at $a$ and $b$ respectively.  This leads to $2H = 0.0666$.  We note again that our algorithm is also somewhat sensitive to the boundary data in this case, so we choose the coefficients carefully.  Also, we used an initial guess of uniform height 0 again as that performed better than one with uniform height 1.

\begin{figure}[h!]
	\centering
	\scalebox{0.195}{\includegraphics{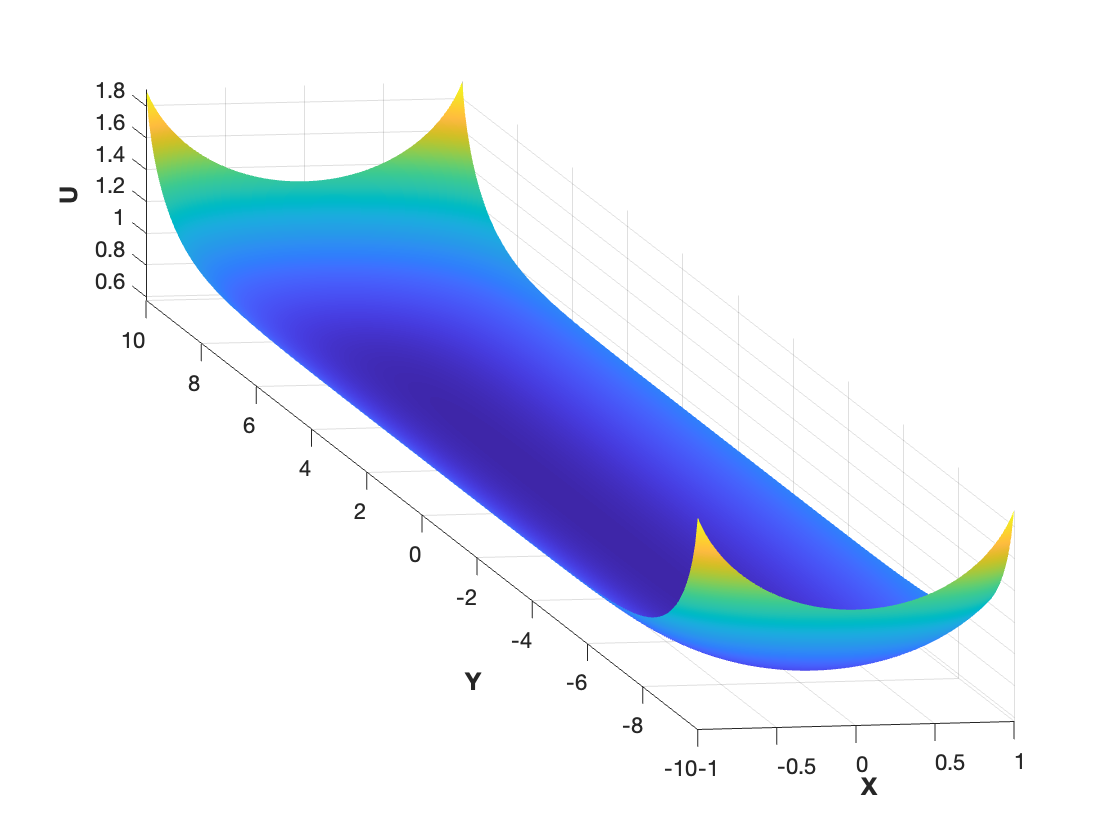}}
	\scalebox{0.18}{\includegraphics{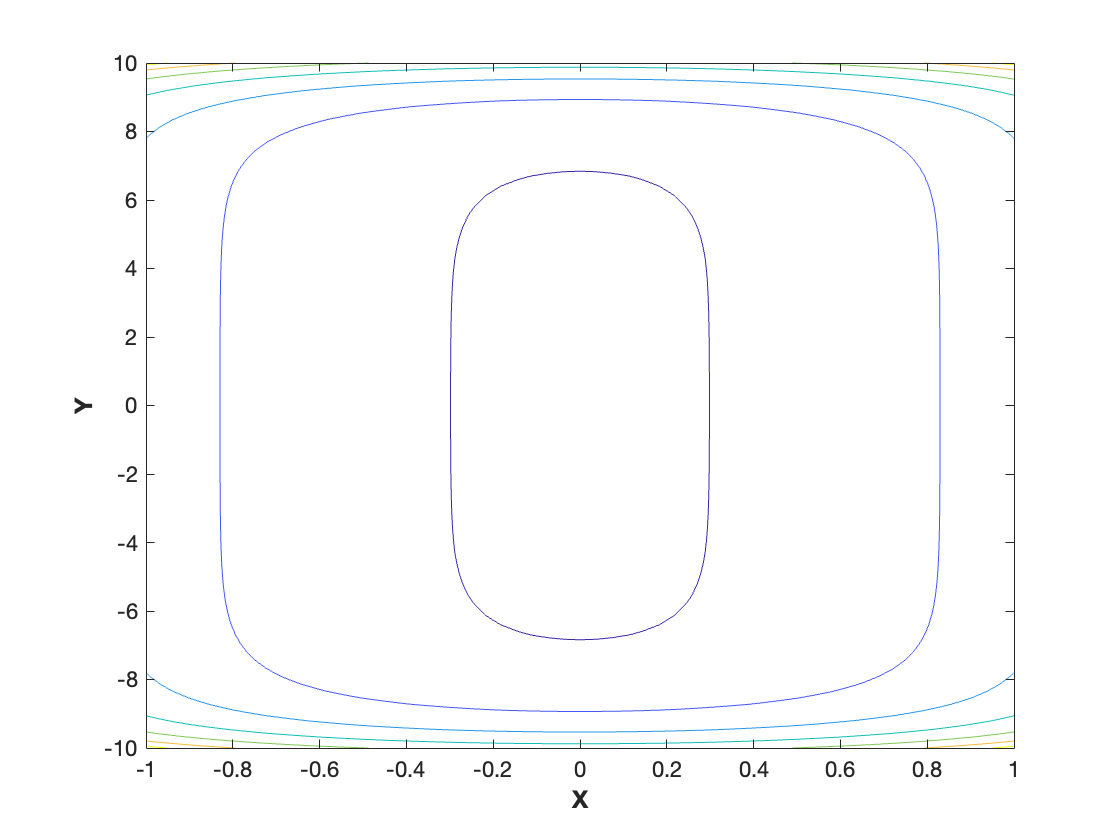}}
	\scalebox{0.195}{\includegraphics{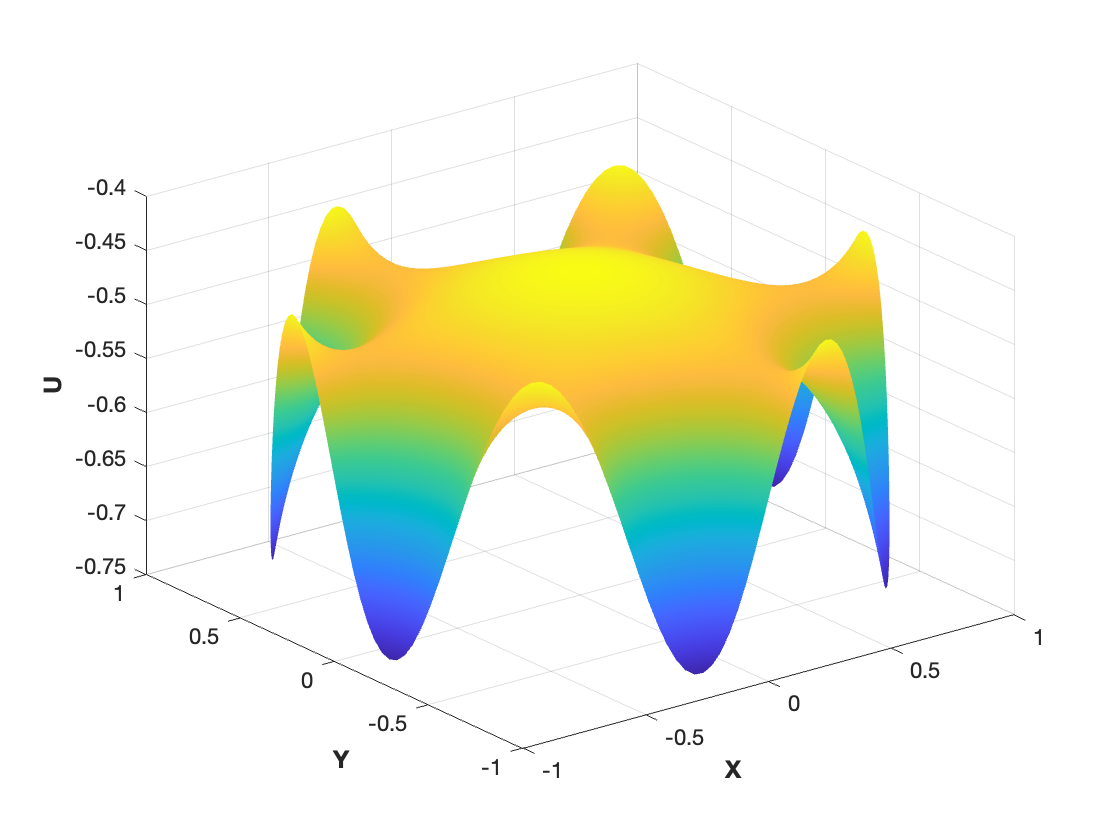}}
	\scalebox{0.195}{\includegraphics{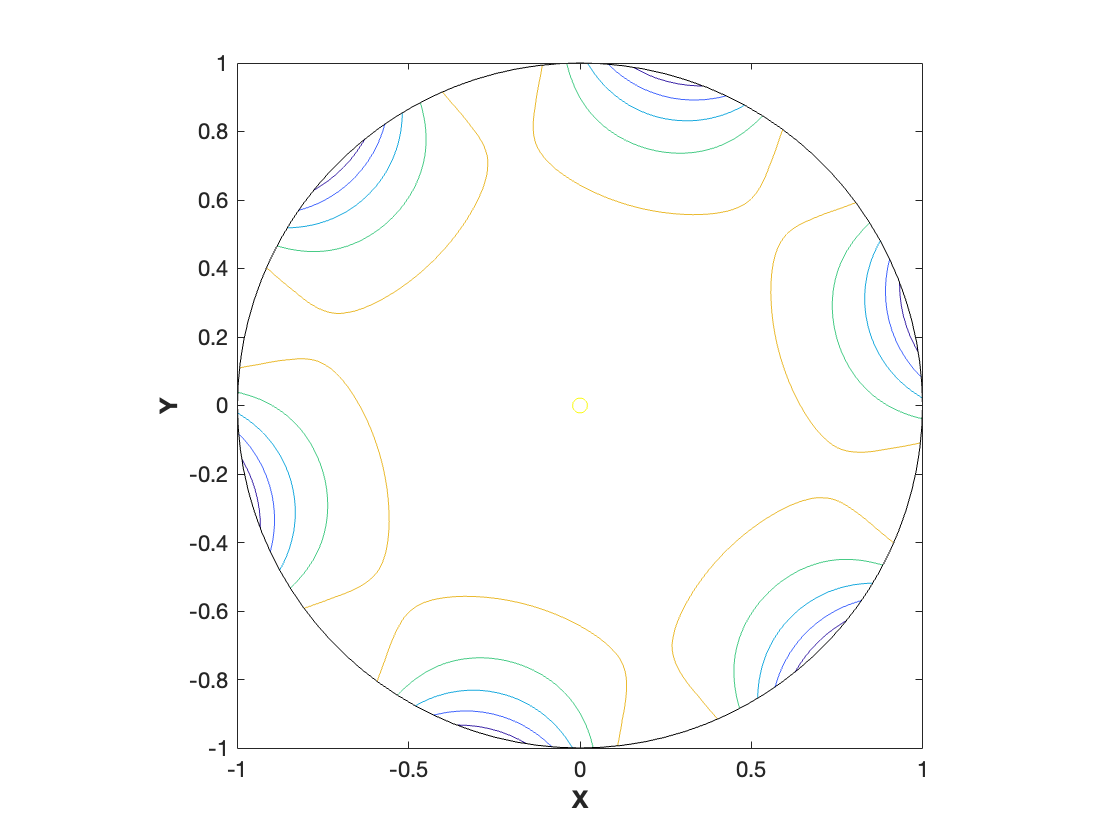}}
	\scalebox{0.195}{\includegraphics{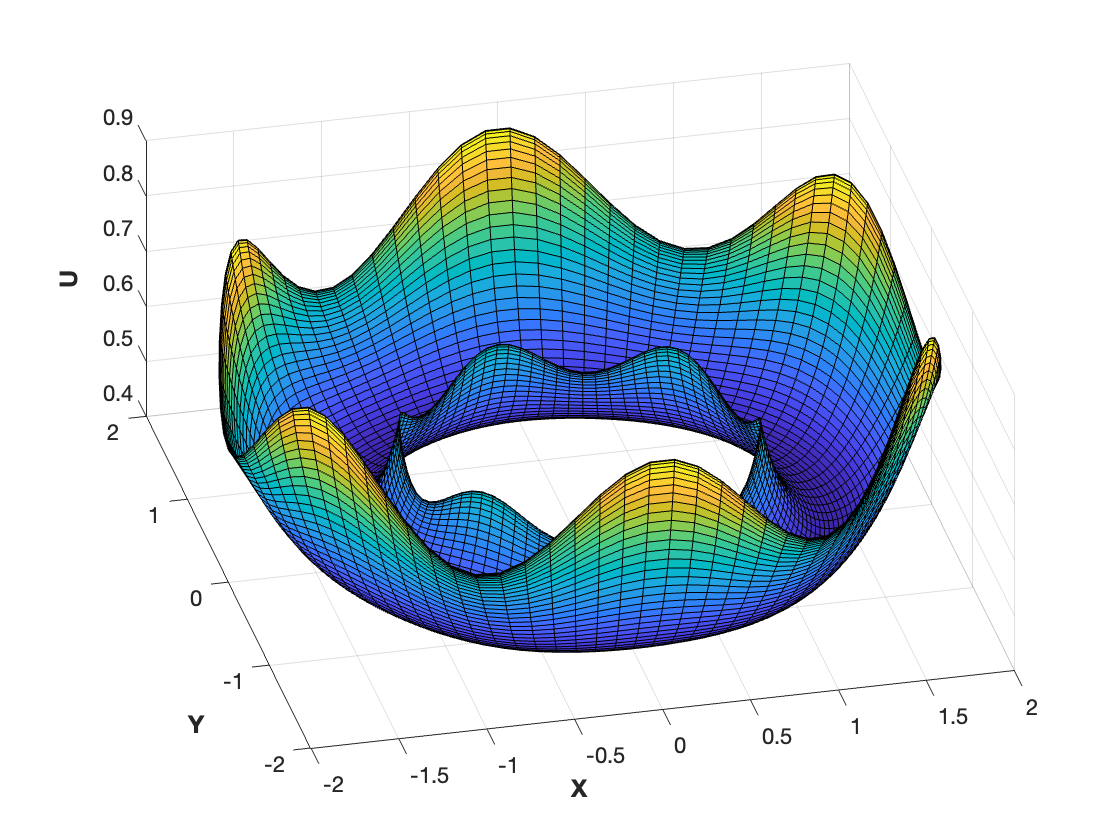}}
	\scalebox{0.195}{\includegraphics{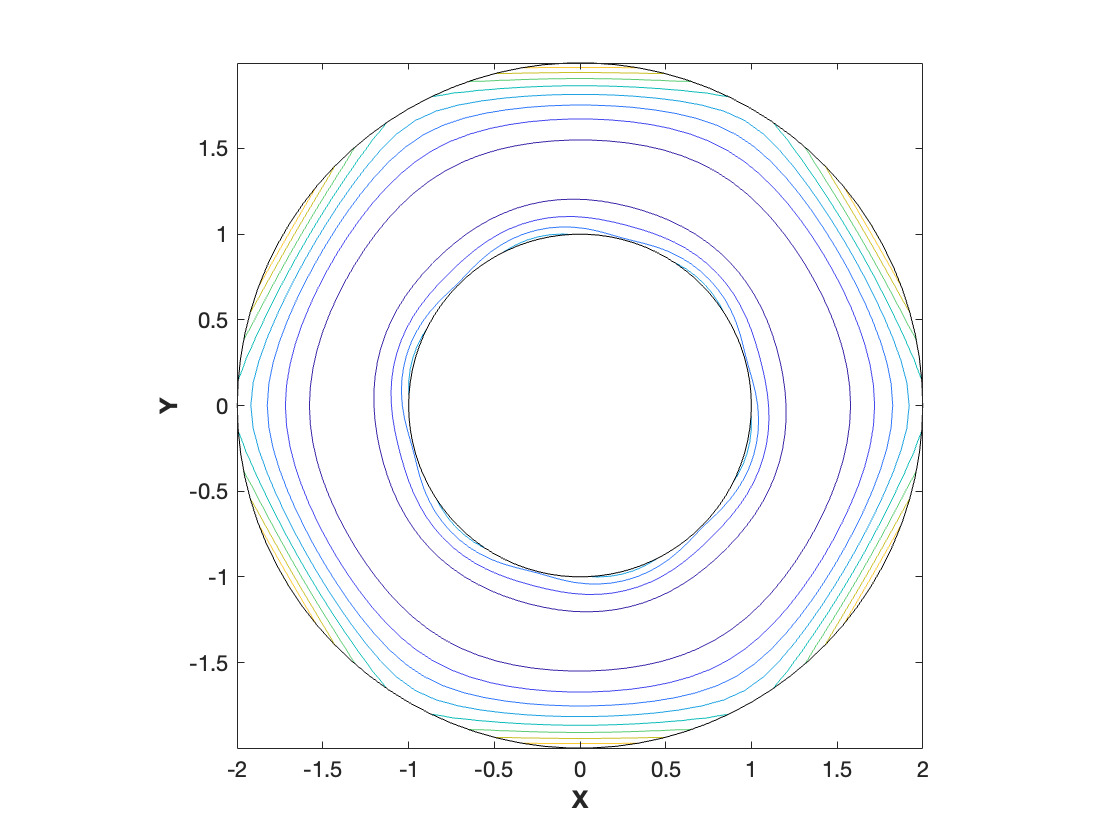}}
	\caption{Capillary surfaces on three different domains.  }
	\label{fig:CapNeu}
\end{figure}

\begin{figure}[h!]
	\centering
	\scalebox{0.195}{\includegraphics{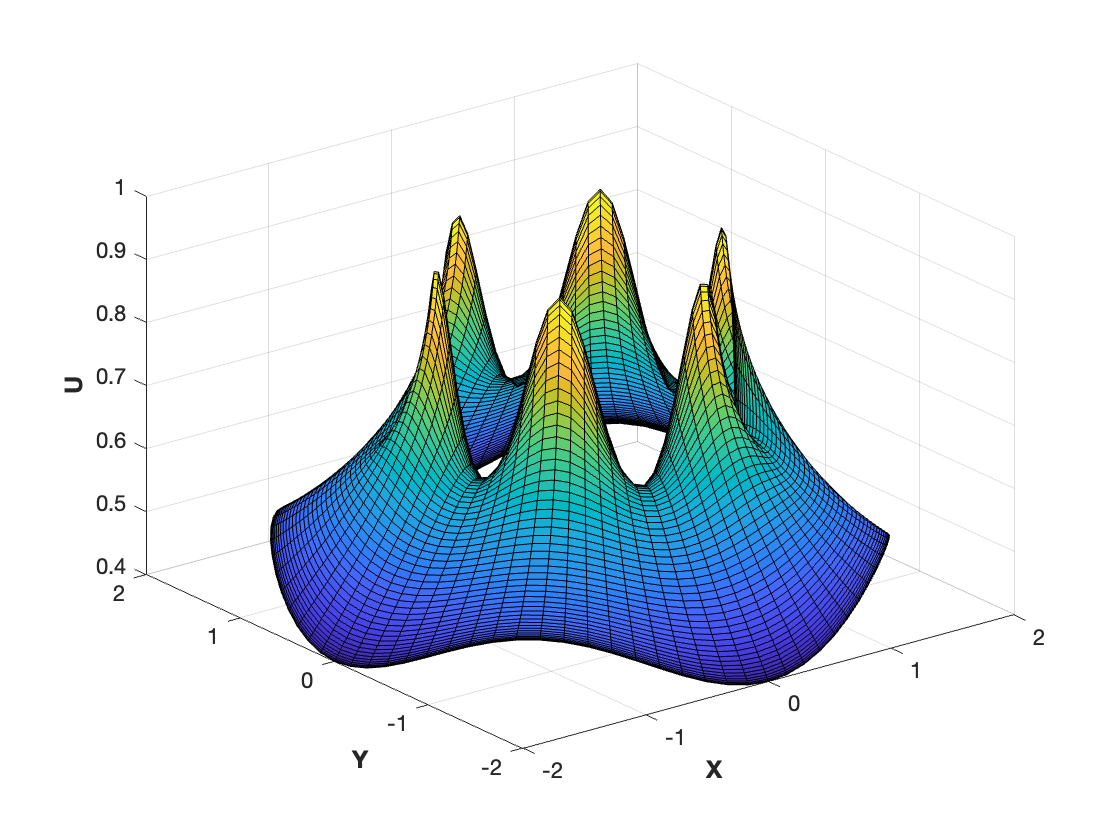}}
	\scalebox{0.195}{\includegraphics{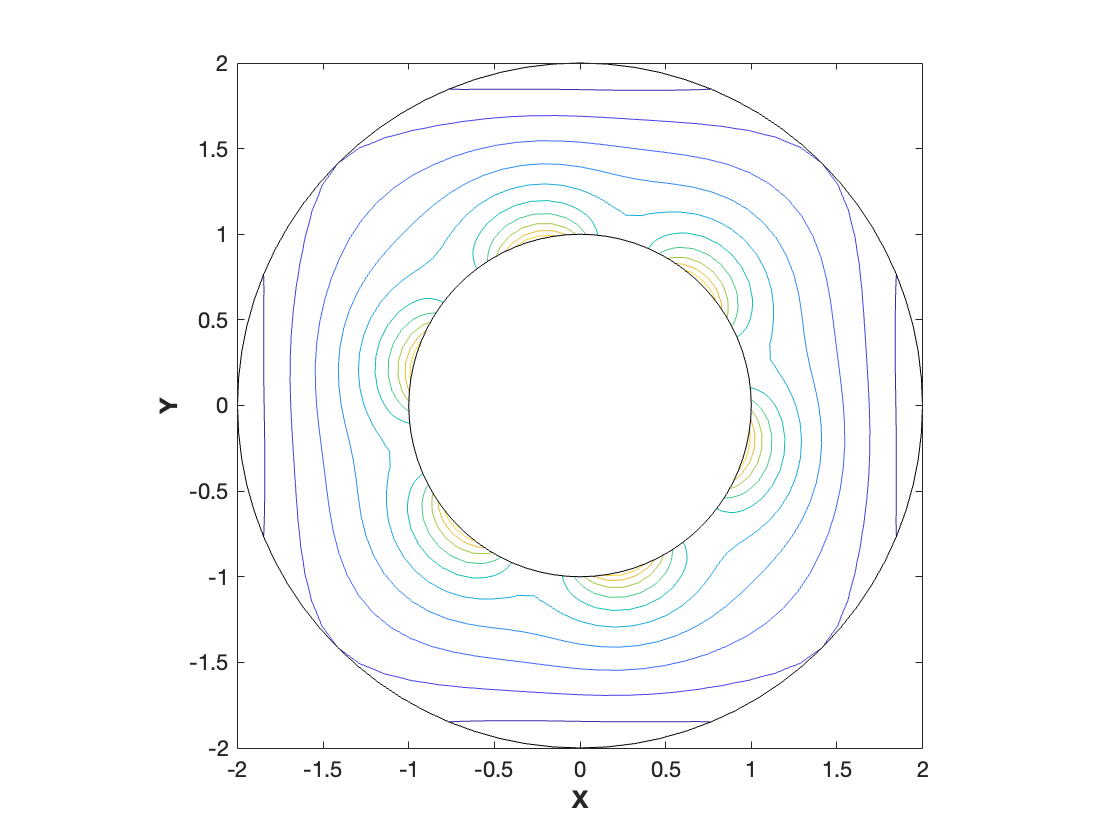}}
	\scalebox{0.195}{\includegraphics{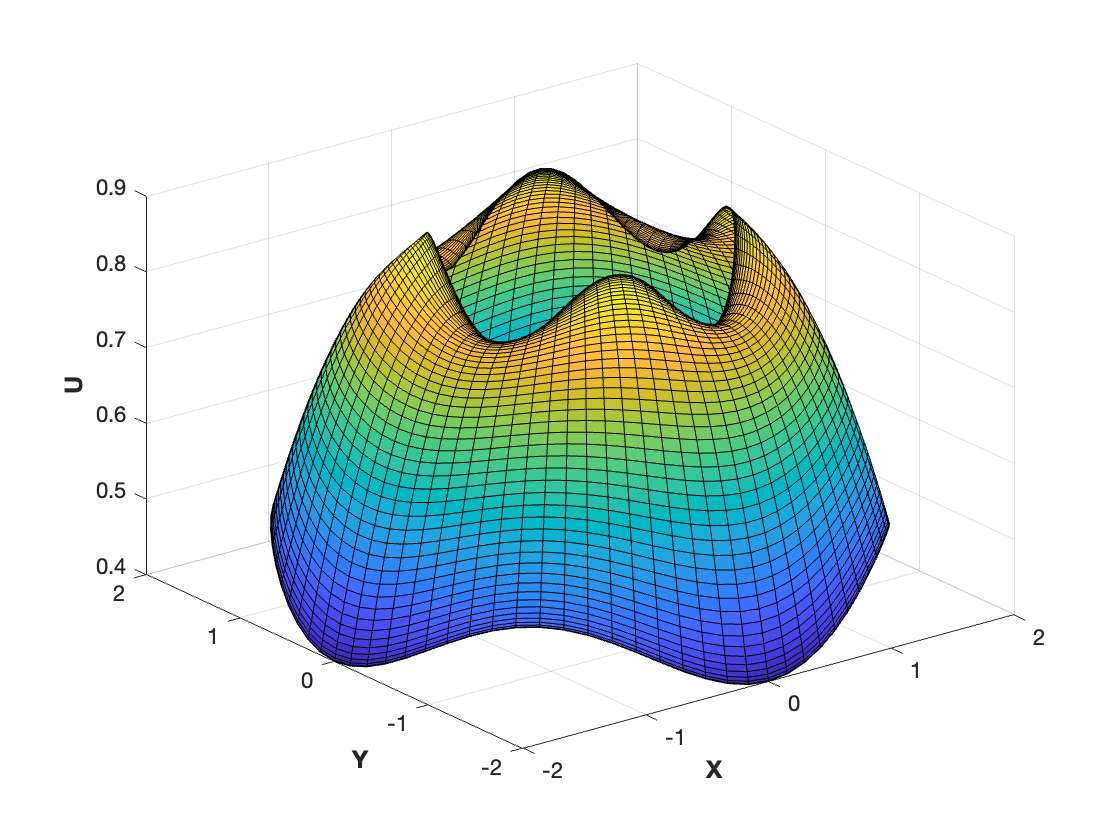}}
	\scalebox{0.195}{\includegraphics{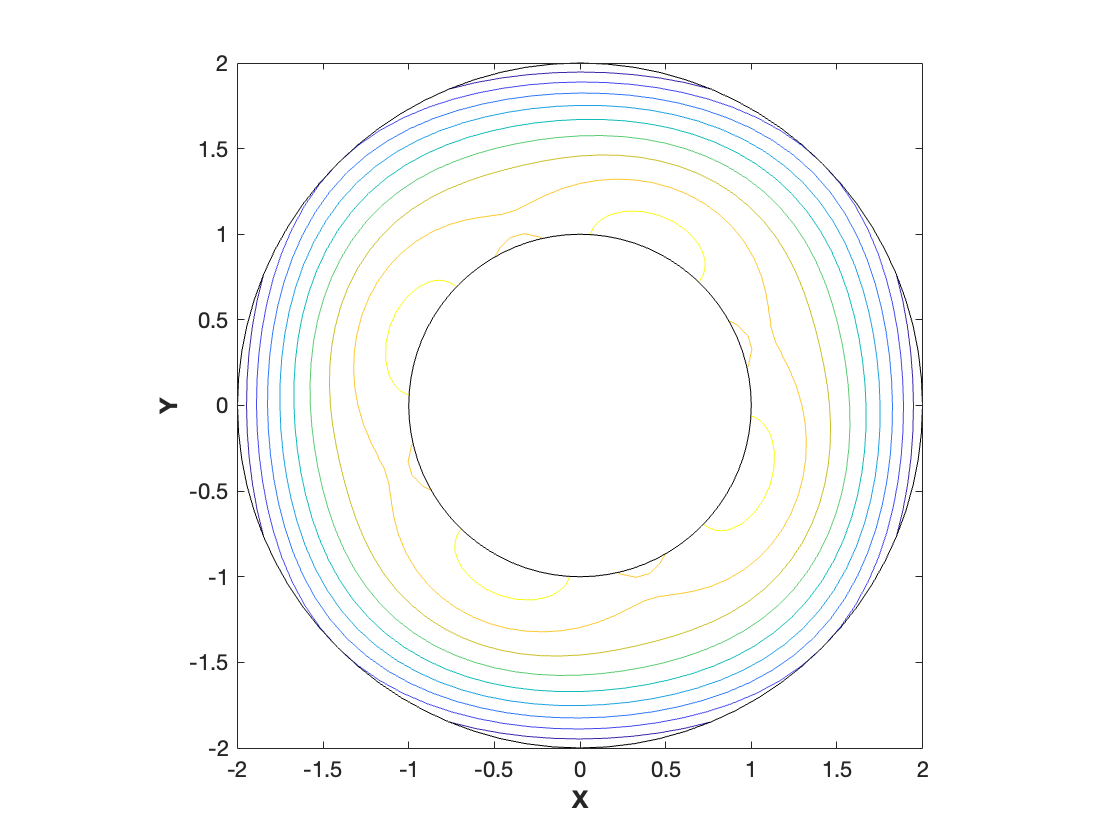}}
	\scalebox{0.195}{\includegraphics{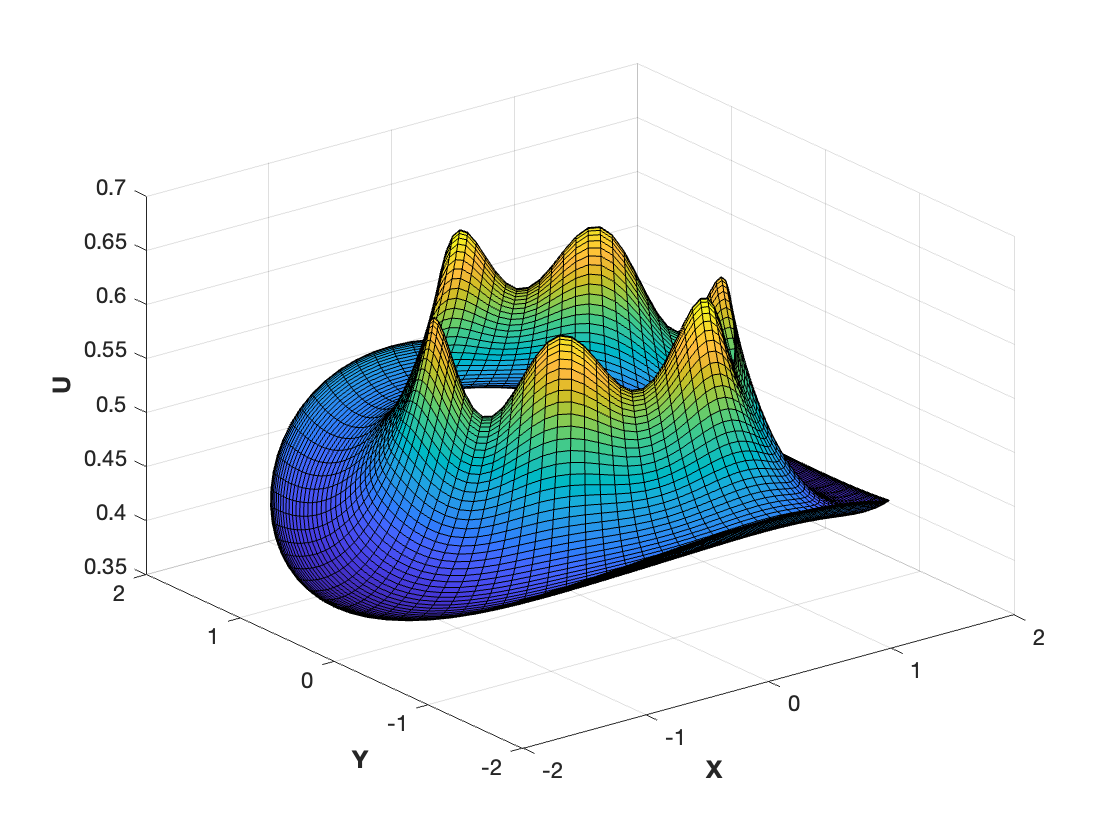}}
	\scalebox{0.195}{\includegraphics{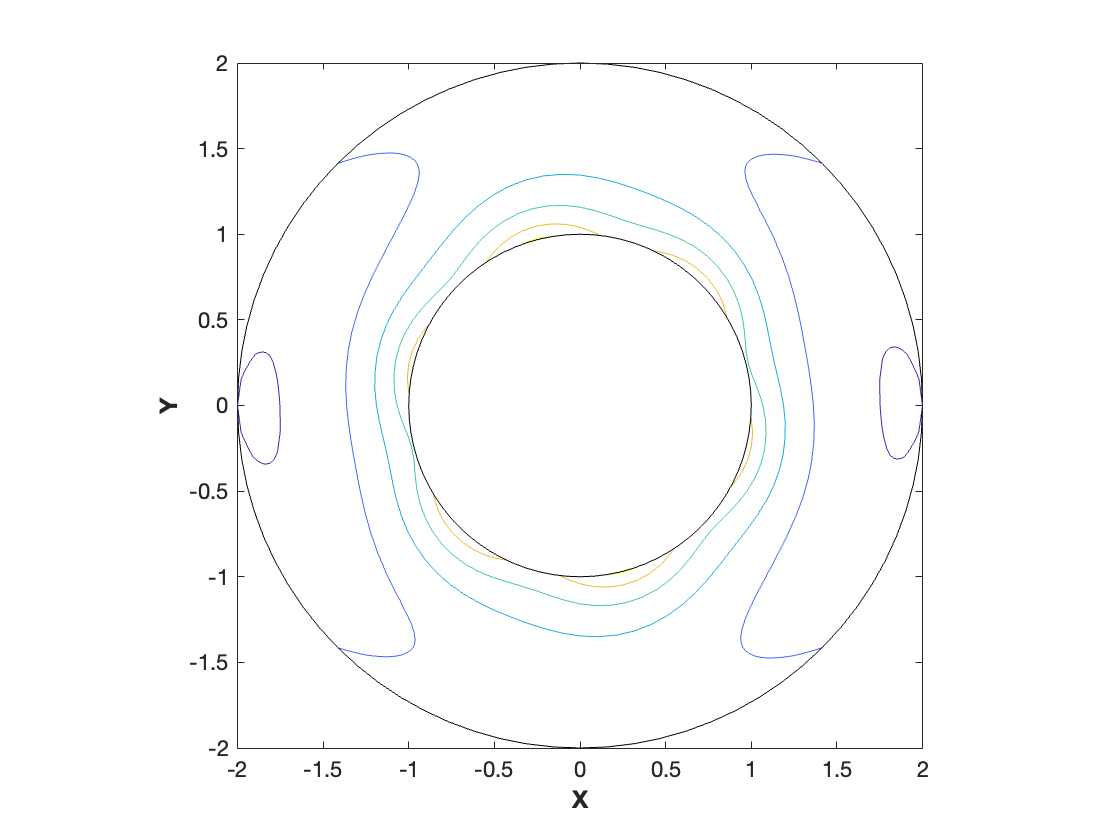}}
	\caption{A minimal surface, a CMC surface, and a capillary surface  are all shown with Dirichlet data on the outer boundary and with capillary data on the inner radius.}
	\label{fig:Mix}
\end{figure}

Figure~\ref{fig:CapNeu} shows examples of the capillary surfaces with capillary data on the boundary for three different domains.  The top figure is a repeat of the data $\gamma = \pi/4 + 0.035$ from Figure~\ref{fig:CapCorner} where here we take the domain to be $[-1,1]\times[-10,10]$.  In the square case earlier in the paper we can use $n=40$ grid points and the convergence takes 8 seconds.  For the non-square case we got $n=87$ grid points after the adaptive algorithm added points, and the computation took 724 seconds, which is more than 12 minutes.  For the middle figure with a disk domain we used the contact angle data of
$$
\gamma(\theta) = \pi/2 + 0.3 + 0.75\sin(6\theta).
$$
For the annular domain in the bottom figure the boundary data was
$$
\gamma_a(\theta) =  \pi/3 + 0.2\sin(6\theta),
$$
and
$$
\gamma_b(\theta) =  \pi/3 + 0.2\cos(6\theta)
$$
on the inner and outer radii of $a = 1$ and $b = 2$ respectively.  

Figure~\ref{fig:Mix} shows examples of a minimal surface, a CMC surface, and a capillary surface all with a mix of Dirichlet and capillary data on the boundaries of an annulus with radii $a = 1$ and $b = 1$.  The top figure shows a minimal surface with capillary data of 
$$
\gamma_a(\theta) =  \pi/3 + 0.75\sin(6\theta),
$$
at $a$ and Dirichlet data of
$$
g_b(\theta) = 1/2 - 0.1\cos^2(2\theta)
$$
at $b$.    The middle figure shows a CMC surface with capillary data of 
$$
\gamma_a(\theta) =  \pi/2 + 0.2\sin(4\theta),
$$
at $a$ and Dirichlet data of
$$
g_b(\theta) = 1/2 - 0.1\cos^2(2\theta)
$$
at $b$ and the mean curvature was $2H = -0.85$.   The bottom figure shows a capillary surface with capillary data of 
$$
\gamma_a(\theta) =  \pi/3 + 0.2\sin(6\theta),
$$
at $a$ and Dirichlet data of
$$
g_b(\theta) = 1/2 - 0.1\cos^2(\theta)
$$
at $b$.



\begin{bibdiv}
		\begin{biblist}
			
										\bib{AurentzTrefethen2017}{article}{
				author={Aurentz, Jared L.},
				author={Trefethen, Lloyd N.},
				title={Block operators and spectral discretizations},
				journal={SIAM Rev.},
				volume={59},
				date={2017},
				number={2},
				pages={423--446},
				issn={0036-1445},
				review={\MR{3646500}},
				doi={10.1137/16M1065975},
			}
			

			
			\bib{BirkissonDriscoll2012}{article}{
				author={Birkisson, Asgeir},
				author={Driscoll, Tobin A.},
				title={Automatic Fr\'{e}chet differentiation for the numerical solution of
					boundary-value problems},
				journal={ACM Trans. Math. Software},
				volume={38},
				date={2012},
				number={4},
				pages={Art. 26, 29},
				issn={0098-3500},
				review={\MR{2972670}},
				doi={10.1145/2331130.2331134},
			}
			
						\bib{Bliss1925}{book}{
				title={Calculus of variations},
				author={Bliss, Gilbert Ames},
				number={1},
				date={1925},
				publisher={Mathematical Association of America}
			}
			
			\bib{Boyd2001}{book}{
				author={Boyd, John P.},
				title={Chebyshev and Fourier spectral methods},
				edition={2},
				publisher={Dover Publications, Inc., Mineola, NY},
				date={2001},
				pages={xvi+668},
				isbn={0-486-41183-4},
				review={\MR{1874071}},
			}
			
			\bib{Brubaker2018}{article}{
				author={Brubaker, Nicholas D.},
				title={A continuation method for computing constant mean curvature
					surfaces with boundary},
				journal={SIAM J. Sci. Comput.},
				volume={40},
				date={2018},
				number={4},
				pages={A2568--A2583},
				issn={1064-8275},
				review={\MR{3845275}},
				doi={10.1137/17M1143228},
			}
			
			\bib{Burns2022}{article}{
			author={Burns, Keaton J. },
			author={Fortunato, Daniel},
			author={Julien, Keith },
			author={Vasil, Geoffrey M. },
			title={Corner Cases of the Generalized Tau Method}
    		journal={preprint arXiv:2211.17259},
    		date={2022},
    	}
			

			%
			
			\bib{Chebfun}{book}{
				editor={Driscoll, T. A.},
				editor={Hale, N.},
				editor={Trefethen, L. N.},
				title={Chebfun Guide},
				publisher={Pafnuty Publications},
				place={Oxford},
				date={2014},
			}
			
			
			
			\bib{Finn1965}{article}{
				author={Finn, Robert},
				title={Remarks relevant to minimal surfaces, and to surfaces of
					prescribed mean curvature},
				journal={J. Analyse Math.},
				volume={14},
				date={1965},
				pages={139--160},
				issn={0021-7670},
				review={\MR{0188909}},
				doi={10.1007/BF02806384},
			}
			
			
			\bib{ecs}{book}{
				author={Finn, Robert},
				title={Equilibrium capillary surfaces},
				series={Grundlehren der Mathematischen Wissenschaften [Fundamental
					Principles of Mathematical Sciences]},
				volume={284},
				publisher={Springer-Verlag},
				place={New York},
				date={1986},
				pages={xvi+245},
				isbn={0-387-96174-7},
			}
			
			\bib{Fornberg1996}{book}{
				author={Fornberg, Bengt},
				title={A practical guide to pseudospectral methods},
				series={Cambridge Monographs on Applied and Computational Mathematics},
				volume={1},
				publisher={Cambridge University Press, Cambridge},
				date={1996},
				pages={x+231},
				isbn={0-521-49582-2},
				review={\MR{1386891}},
				doi={10.1017/CBO9780511626357},
			}
			
			\bib{Gauss}{book}{
				author={C. F. Gauss},
				title={Principia generalia theoriae figurae fluidorum in statu aequilibrii},
				publisher={Dieterich, Göttingen, Germany},
				date={1830}
			}
			
				\bib{Giusti1984}{book}{
				author={Giusti, Enrico},
				title={Minimal surfaces and functions of bounded variation},
				series={Monographs in Mathematics},
				volume={80},
				publisher={Birkh\"{a}user Verlag, Basel},
				date={1984},
				pages={xii+240},
				isbn={0-8176-3153-4},
				review={\MR{0775682}},
				doi={10.1007/978-1-4684-9486-0},
			}
			
			\bib{GraggTapia1974}{article}{
				author={Gragg, W. B.},
				author={Tapia, R. A.},
				title={Optimal error bounds for the Newton-Kantorovich theorem},
				journal={SIAM J. Numer. Anal.},
				volume={11},
				date={1974},
				pages={10--13},
				issn={0036-1429},
				review={\MR{0343594}},
				doi={10.1137/0711002},
			}
			
			\bib{HaugTreinen2024}{article}{
				author={Haug, Jonas},
				author={Treinen, Ray},
				title={Multi-scale spectral methods for bounded radially symmetric
					capillary surfaces},
				journal={Electron. Trans. Numer. Anal.},
				volume={60},
				date={2024},
				pages={20--39},
				review={\MR{4695961}},
				doi={10.1553/etna\_vol60s20},
			}
			
			\bib{KantorovichAkilov1982}{book}{
				author={Kantorovich, L. V.},
				author={Akilov, G. P.},
				title={Functional analysis},
				edition={2},
				note={Translated from the Russian by Howard L. Silcock},
				publisher={Pergamon Press, Oxford-Elmsford, N.Y.},
				date={1982},
				pages={xiv+589},
				isbn={0-08-023036-9},
				review={\MR{0664597}},
			}
			
			\bib{Lopez2013}{book}{
				author={L\'{o}pez, Rafael},
				title={Constant mean curvature surfaces with boundary},
				series={Springer Monographs in Mathematics},
				publisher={Springer, Heidelberg},
				date={2013},
				pages={xiv+292},
				isbn={978-3-642-39625-0},
				isbn={978-3-642-39626-7},
				review={\MR{3098467}},
				doi={10.1007/978-3-642-39626-7},
			}
%
%
%
%

\bib{NocedalWright2006}{book}{
	author={Nocedal, Jorge},
	author={Wright, Stephen J.},
	title={Numerical optimization},
	series={Springer Series in Operations Research and Financial Engineering},
	edition={2},
	publisher={Springer, New York},
	date={2006},
	pages={xxii+664},
	isbn={978-0387-30303-1},
	isbn={0-387-30303-0},
	review={\MR{2244940}},
}

\bib{Ortega1968}{article}{
	author={Ortega, James M.},
	title={The Newton-Kantorovich theorem},
	journal={Amer. Math. Monthly},
	volume={75},
	date={1968},
	pages={658--660},
	issn={0002-9890},
	review={\MR{0231218}},
	doi={10.2307/2313800},
}
			
			\bib{Ros1996}{article}{
				author={Ros, Antonio},
				title={Embedded minimal surfaces: forces, topology and symmetries},
				journal={Calc. Var. Partial Differential Equations},
				volume={4},
				date={1996},
				number={5},
				pages={469--496},
				issn={0944-2669},
				review={\MR{1402733}},
				doi={10.1007/s005260050050},
			}
			
			\bib{Serrin1969}{article}{
				author={Serrin, J.},
				title={The problem of Dirichlet for quasilinear elliptic differential
					equations with many independent variables},
				journal={Philos. Trans. Roy. Soc. London Ser. A},
				volume={264},
				date={1969},
				pages={413--496},
				issn={0080-4614},
				review={\MR{0282058}},
				doi={10.1098/rsta.1969.0033},
			}
			
			
			\bib{Trefethen2000}{book}{
				author={Trefethen, Lloyd N.},
				title={Spectral methods in MATLAB},
				series={Software, Environments, and Tools},
				volume={10},
				publisher={Society for Industrial and Applied Mathematics (SIAM),
					Philadelphia, PA},
				date={2000},
				pages={xviii+165},
				isbn={0-89871-465-6},
				review={\MR{1776072}},
				doi={10.1137/1.9780898719598},
			}
			
			
							\bib{Treinen2023a}{article}{
									author={Treinen, Ray},
									title={Spectral methods for capillary surfaces described by bounded
											generating curves},
									journal={Appl. Math. Comput.},
									volume={450},
									date={2023},
									pages={Paper No. 127886, 17},
									issn={0096-3003},
									review={\MR{4566044}},
									doi={10.1016/j.amc.2023.127886},
								}
								
								\bib{Treinen2023b}{article}{
									author={Treinen, Raymond},
									title={Discussion of a uniqueness result in ``equilibrium configurations
										for a floating drop''},
									journal={Electron. J. Differential Equations},
									date={2023},
									pages={Paper No. 32, 11},
									review={\MR{4574300}},
									doi={10.58997/ejde.2023.32},
								}
			
			\bib{WilberTownsendWright2017}{article}{
				author={Wilber, Heather},
				author={Townsend, Alex},
				author={Wright, Grady B.},
				title={Computing with functions in spherical and polar geometries II. The
					disk},
				journal={SIAM J. Sci. Comput.},
				volume={39},
				date={2017},
				number={3},
				pages={C238--C262},
				issn={1064-8275},
				review={\MR{3666775}},
				doi={10.1137/16M1070207},
			}
			
			\bib{Zeidler1986}{book}{
				author={Zeidler, Eberhard},
				title={Nonlinear functional analysis and its applications. I},
				note={Fixed-point theorems;
					Translated from the German by Peter R. Wadsack},
				publisher={Springer-Verlag, New York},
				date={1986},
				pages={xxi+897},
				isbn={0-387-90914-1},
				review={\MR{0816732}},
				doi={10.1007/978-1-4612-4838-5},
			}
			
		\end{biblist}
	\end{bibdiv}
\end{document}